\documentclass[a4paper, 12pt]{article}
\usepackage{enumerate, theorem}
\usepackage{amsmath, amsfonts, amssymb}
\usepackage[height=22.5cm, width=15cm]{geometry}
\usepackage[all]{xy}
\usepackage{mathrsfs, yfonts}

\DeclareMathOperator{\C}{\mathbb{C}}
\DeclareMathOperator{\A}{\mathcal{A}_S}

\newcommand{\parag}[1]{\paragraph{\sc{#1.}}}

\newtheorem{thm}{Theorem}[subsection]
\newtheorem{defn}[thm]{Definition}
\newtheorem{cor}[thm]{Corollary}
\newtheorem{prop}[thm]{Proposition}
\newtheorem{lemma}[thm]{Lemma}

\begin{document}

\title{A  finiteness theorem for \ $S-$relative formal  Brieskorn modules.\\ version 2}

\author{Daniel Barlet\footnote{Barlet Daniel, Institut Elie Cartan UMR 7502  \newline
Universit\'e de Lorraine, CNRS, INRIA  et  Institut Universitaire de France, \newline
BP 239 - F - 54506 Vandoeuvre-l\`es-Nancy Cedex.France. \newline
e-mail : daniel.barlet@univ-lorraine.fr}.}

\date{5/02/14.}

\maketitle

\section*{Abstract}

We give a general result of finiteness for holomorphic families of Brieskorn modules constructed from a holomorphic family of one parameter degeneration of compact complex manifolds acquiring (general) singularities.

\parag{AMS Classification} 32 S 25, 32 S 40, 32 S 50.

\parag{Key words} Brieskorn module,  (a,b)-module, asymptotic expansion, Gauss-Manin connection, filtered differential equation.

\section*{Introduction}

In this article we are interested in the following situation :\\
Let \ $f : X \to T$ \ be a proper  holomorphic and surjective map between complex manifolds such that outside an hypersurface \ $S \subset T$ \ the fibers of \ $f$ \ are smooth. Then around the generic point of \ $S$ \ we may assume that \ $T $ \ is locally isomorphic to \ $S_0\times D$ \ where \ $D$ \ is a small open disc with center \ $0$ \ in \ $\C$ \ and \ $S_0$ \ an open set in \ $S$. We can consider the restriction of 
\ $f$ \ over \ $S_0\times D$ \ as a holomorphic  family parametrized by \ $S_0$ \ of one parameter degeneration of compact complex manifolds acquiring singularities over \ $S_0\times \{0\}$. \\
Our  first  result is to give in such a local  situation (with suitable hypotheses, always satisfied in the absolute case and in generic relative situations) a \ $S-$relative version of the construction given in [B. II] theorem 2.1.1. Recall that this result gives a complex of \ $\mathcal{A}-$modules quasi-isomorphic to the formal completion in \ $f$ \ of the complex \ $(Ker\,df^{\bullet}, d^{\bullet})$. The action of the algebra \ $\mathcal{A}$ \  realizing the filtered Gauss-Manin connection of \ $f$. We obtain the relative case in theorem \ref{construction}.\\

 In the second part we show that the properness assumption allows to obtain a finiteness result for the family parametrized by \ $S$ \ of  Brieskorn modules (see [Br.70], [S.89], [B.93] and [B.95]) obtained via the cohomology of the direct image of the complex constructed in the local setting. This is proved in the theorem \ref{Finitude}  \\
The finiteness  result, which is new already in the absolute case, will be a key tool to produce holomorphic families of Brieskorn modules when we consider holomorphic families of one parameter  degenerations of compact complex manifolds acquiring arbitrary singularities.\\
Note also that our hypotheses for the theorem \ref{construction} are purely local and that we use properness without any K\"ahler type assumption for the finiteness theorem \ref{Finitude}.

\tableofcontents

\section{ The general  construction of (a,b)-modules in the relative case. }

\subsection{Our situation.}

 \parag{Notations} Let \ $S$ \ be a reduced complex space. We shall say that \ $\pi : \mathcal{X} \to S$ \ is a \ {\bf $S-$manifold} when \ $\mathcal{X}$ \ is a reduced complex space and \ $\pi$ \ a holomorphic map which is \ $S-$smooth. By definition this means that locally on \ $\mathcal{X}$ \ we have a \ $S-$isomorphism of \ $\mathcal{X}$ \ with a product \ $S \times U$ \ where \ $U$ \ is an open set in  \ $\C^{n+1}$. Such an isomorphism will be called a \ $S-$relative chart or a \ $S-$relative system of coordinates on \ $\mathcal{X}$. \\
 
 In such a situation we have on \ $\mathcal{X}$ \ a locally free \ $\mathcal{O}_{\mathcal{X}}-$sheaf of \ {\bf $S-$relative holomorphic differential forms} \ $\Omega_{/S}^{\bullet}$ \ which corresponds, via the local \ $S-$isomorphisms above, to the sheaf of holomorphic differential forms on \ $U$ \ with holomorphic coefficients  on \ $S \times U$.\\
 For a holomorphic function \ $f$ \ on \ $\mathcal{X}$ \ the \ {\bf $S-$relative differential} is defined as a section of the sheaf \ $\Omega^1_{/S}$ \  and denoted \ $d_{/S}f$. If \ $d_{/S}f_x \not= 0$, then \ $f$ \ may be choosen as the first coordinate of a  local \ $S-$relative coordinate system on \ $\mathcal{X}$ \ near the point \ $x$. The relative differential \ $d_{/S}$ \ is well defined  and of degree \ $1$ \ on the graded sheaf \ $\Omega_{/S}^{\bullet}$. The relative differential  \ $d_{/S}$ \ is \ $f^{-1}(\mathcal{O}_S)-$linear and  the \ $S-$relative de Rham  complex is a complex of \ $f^{-1}(\mathcal{O}_S)$ \ modules. The de Rham lemma with holomorphic  parameters gives that the complex \ $(\Omega_{/S}^{\bullet}. d^{\bullet}_{/S})$ \ has zero cohomology sheaves in positive degrees ; so it   is a resolution of \ $f^{-1}(\mathcal{O}_S)$. \\

\parag{The situation  $(@)$} We shall consider a holomorphic function \ $f : \mathcal{X} \to D$ \ where \ $D$ \ is a disc with center \ $0$ \ in \ $\C$ \ and where \ $\pi : \mathcal{X} \to S$ \ is a \ $S-$manifold.\\
  We shall define \ $\mathcal{Y} : = f^{-1}(0)$ \ and \ $\mathcal{Z} : = \{ x \in \mathcal{X} \ / \ d_{/S}f_x = 0\}$. 
 \parag{Hypothesis \ $(H0)$} We shall assume that \ $\mathcal{Y}$ \ has codimension 1 in each  fiber of \ $\pi$, that \ $ \mathcal{Z} \subset \mathcal{Y}$ \ and that \ $\mathcal{Z}$ \ has codimension \ $\geq 2$ \ in each fiber of \ $\pi$.
 \parag{Hypothesis \ $(H1)$} We shall assume that locally on \ $\mathcal{X}$ \ there exists a proper \ $S-$modification \ $\tau : \tilde{\mathcal{X}} \to \mathcal{X}$ \ such that \ $\tilde{\mathcal{X}}$ \ is \ $S-$smooth and  \ $\tilde{\mathcal{Y}} : = \tau^{-1}(\mathcal{Y})$ \ is locally on \ $\tilde{\mathcal{X}}$ \ a \ $S-$relative  normal crossing divisor. That is to say that \ $f$ \ admits a simultaneous desingularisation over \ $S$.
 \parag{Hypothesis \ $(H2)$} We shall assume that locally on \ $\mathcal{Z}$ \ we shall find local \ $S-$relative Milnor fibration.\\

 The precise meaning of these  hypotheses is given below, but  remark already  that in the absolute case ( i. e. for \ $S = \{0\}$), the hypotheses \ $(H1)$ \ and \ $(H2)$  are always satisfied for any non zero germ of holomorphic function in \ $\C^{n+1}$. The hypothesis \ $(H0)$ \ in this absolute case asks simply that the hypersurface \ $f^{-1}(0)$ \ is reduced. Note also that these hypotheses  are purely local on \ $\mathcal{Y}$.\\
  
Let me recall the following standard definitions which makes precise the hypothesis \ $(H1)$.

\begin{defn} We shall say that the divisor \ $  \mathcal{Y} : = \{ f = 0 \}$ \ in the \ $S-$manifold  \ $\mathcal{X}$ \ is  a {\bf \ $S-$normal crossing divisor} if near each point \ $y \in \mathcal{Y}$ \ there exists a system of local \ $S-$coordinate in \ $\mathcal{X}$ \ such that we have \ $f(s,z) = z^{\alpha}$ \ where \ $\alpha$ \ is in \ $\mathbb{N}^{n+1}$.
\end{defn}

\begin{defn} We shall say that the holomorphic  function \ $f : \mathcal{X} \to D$ \ on the \ $S-$manifold  \ $\mathcal{X}$ \ admits a {\bf simultaneous desingularisation over \ $S$} \ if there exists a proper \ $S-$modification \ $\tau : \tilde{\mathcal{X}} \to \mathcal{X}$ \ where \ $\tilde{\mathcal{X}}$ \ is a \ $S-$manifold and such that \ $\tilde{\mathcal{Y}} : = \{ f\circ \tau = 0 \}$ \ is a \ $S-$normal crossing divisor in \ $\tilde{\mathcal{X}}$.
\end{defn}

Note that by a \ $S-$proper modification we mean that \ $\tau$ \  is a  proper holomorphic map  such that there exists  a   closed analytic subset \ $T  \subset \mathcal{X}$ \ with the following properties :
\begin{enumerate}[i)]
\item For each  $s \in S$ \ the fiber \ $T_s$ \ has empty interior in \ $\mathcal{X}_s : = \pi^{-1}(s)$.
\item The map \ $\tau$ \ induces an isomorphism of \ $\tilde{\mathcal{X}} \setminus \tau^{-1}(T)$ \ to \ $\mathcal{X} \setminus T$.
\end{enumerate}

The definition below, which makes precise the hypothesis \ $(H2)$ \ is may-be less standard.

\begin{defn}\label{Milnor relatif}
In the situation \ $(@)$ \ we say that we have locally on \ $\mathcal{Y}$ \  a Milnor fibration for \ $f$ \ is for each point \ $y \in \mathcal{Y}$ \ there exists an open neighbourhood \ $S'$ \ of \ $\pi(y)$ \ is \ $S$\ and a \ $S'-$relative chart \ $\varphi : X \simeq S'\times U$ \ where \ $U$ \ is an open neighbourhood of the origin in \ $\C^{n+1}$, with the following property:
\begin{itemize}
\item For any \ $\varepsilon > 0$ \ sufficently small, there exists  \ $\eta > 0$ \ such that the restriction of \ $\pi\times f$ \ to \ $X \cap f^{-1}(D_{\eta}) \cap \varphi^{-1}(S' \times B_{\varepsilon}) $ \ induces a \ $S'-$relative \ $\mathscr{C}^N-$fibration for \ $N \gg n$ \ from the complement of \ $\mathcal{Y}$ \ to \ $S'\times D^*_{\eta}$ \ with fiber a manifold \ $F_y$ \ with finite dimensional complex homology. Moreover we ask that this \ $\mathscr{C}^N-$fibration is independant of the choices of \ $\varepsilon \ll 1$ \ and \ $\eta \ll \varepsilon$.
\end{itemize}
\end{defn}

Of course, for \ $S = \{0\}$ \ this is always satisfied thanks to Milnor [Mi.68]. In general this is a quite strong condition on the situation \ $(@)$ \ which implies, for instance, the fact that the \ $S-$relative cohomology of the fibers of \ $f$ \ define a local system on \ $S \times D^*$ \ near \ $\pi(y)\times f(y)$. Theen the corresponding  relative Gauss-Manin meromorphic connection has along \ $S\times \{0\}$ \ a regular singularity thanks to [D. 70].

\parag{Question} Does the condition \ $(H1)$ \ implies the condition \ $(H2)$ ?

\subsection{The sheaf \ $\mathcal{A}_S$ \ on a reduced complex space.}

In this paragraph 1.2  we consider the situation \ $(@)$ \ with  the hypothesis \ $(H0)$ \ alone.\\

Let \ $S$ \ be a reduced complex space. We define the following sheaves on \ $S$:
\begin{enumerate}
\item  The sheaf  \ $\mathcal{O}_S[[a]]$ \ of commutative \ $\mathcal{O}_S-$algebras defined as the presheaf
$$ U \mapsto \Gamma(U, \mathcal{O}_S)[[a]] : = \prod_{j=0}^{\infty} \ \Gamma(U,\mathcal{O}_S).a^j .$$
We shall denote \ $\big(\mathcal{O}_S[[a]]\big)_{s_0}$ \ the germ at a point \ $s_0 \in S$ \ to distinguish it from the algebra \ $\mathcal{O}_{S,s_0}[[a]]$ \ which is strictely bigger in general\footnote{For an element in \ $\big(\mathcal{O}_S[[a]]\big)_{s_0}$ \ the coefficients of the formal power serie in \ $a$ \ are defined and holomorphic on a common open neighbourghood of \ $s_{0}$ \ in \ $S$.} . Note that both are noetherain rings but the proof of this fact for the first one uses Oka's theorem to have the coherence of \ $\mathcal{O}_S$ \ and Cartan theorem B on Stein neighbourhoods of the point \ $s_0$ \ in \ $S$. 
\item  The sheaf \ $\mathcal{A}_S^0 $ \  of non commutative \ $\mathcal{O}_S-$algebras defined as the presheaf
$$ U \mapsto  \Gamma(U, \mathcal{O}_S)[a]< b>  : = \oplus_{j = 0}^{\infty} \ .b^j.\Gamma(U, \mathcal{O}_S)[a]$$
with the commutation relation \ $a.b - b.a = b^2$ \ and more generally with the relation  \ $T(a).b = b.T(a) + b.T'(a).b $ \ where \ $T'(a)$ \ is the usual derivative in \ $a$ \ of the  polynomial \ $T(a) \in \Gamma(U, \mathcal{O}_S)[a]$.
\item The sheaf \ $\mathcal{A}_S $ \   of non commutative \ $\mathcal{O}_S-$algebras defined as the presheaf
$$ U \mapsto  \Gamma(U, \mathcal{O}_S)[[a]]<< b>>  : = \prod_{j = 0}^{\infty} \ b^j.\Gamma(U, \mathcal{O}_S)[[a]] $$
with the commutation relation \ $T(a).b = b.T(a) + b.T'(a).b $ \ and the fact that the right and left actions of \ $\mathcal{O}_S[[a]]$ \ are continuous with respect to the \ $b-$adic filtration. 
\end{enumerate}

 \parag{Remark} The sheaf \ $\mathcal{A}_S $ \ is the formal completion (in \ $a$, that is to say in a local coordinate near \ $0$ \  in \ $D$) \ of the sheaf of \ $S-$relative formal microdifferential operators of order \ $\leq 0$ \ on \ $S \times D$. It has a natural action on the sheaf \ $\mathcal{O}_S[[a]]$ \ which is given by the action of \ $b$ \ defined as  \ $[b.g](s,a) : = \int_0^a \ g(s,x).dx $. It is easy to see that it extends as an \ $\mathcal{O}_S[[a]]-$linear action (see the lemma 1.2.3).
 
   \parag{Notations} Let \ $E$ \ be a left \ $\mathcal{A}^0_S-$module. We define \ $A(E)$ \ and \ $B(E)$ \ respectively as the $a-$torsion and \ $b-$torsion of \ $E$. Remark that \ $B(E)$ \ is a \  $\mathcal{A}^0_S-$submodule of \ $E$ \ but that \ $A(E)$ \ is not stable by \ $b$ \ in general. So we shall  define also \ $\tilde{A}(E)$ \ as the maximal \ $\mathcal{A}_S^0-$submodule contained in \ $A(E)$. $\hfill \square$\\
   
   \parag{The \ $\mathcal{A}_S-$module of asymptotic expansions} Let \ $\Lambda$ \ be a finite subset in \ $]0,1] \cap \mathbb{Q}$, and \  $k \geq 0$ \ an integer. Define \ $ \Xi_{\Lambda,S}^{(k)}$ \ as the free \ $\mathcal{O}_S[[b]]-$module generated by \ $e_0(\lambda), \dots,e_k(\lambda)$, where \ $ \lambda$ \ is in \ $ \Lambda$, with the action of \ $a$ \ given by
   $$ a.e_j(\lambda) =  \lambda.b.e_j(\lambda) + b.e_{j-1}(\lambda) $$
   with the convention that \ $e_{-1}(\lambda) = 0 $ \ for each \ $\lambda \in \Lambda$.\\
   If we think of 
    \ $e_j(\lambda)$ \ as 
     $$x^{\lambda-1}.\frac{(Log\, x)^j}{j!} $$
     with \ $b$ \ acting as \ $ \int_0^x $ \  and \ $a$ \ acting as \ $  \times x$ \ it is easy to see that \ $\Xi_{\Lambda,S}^{(k)}$ \ is the standard sheaf with holomorphic coefficients in \ $\mathcal{O}_S$ \ for multivalued  asymptotic expansions in \ $x$ \ with logarithmic terms of degree at most \ $k$,  monodromy with spectrum in \ $\exp(2i\pi.\lambda), \lambda \in \Lambda$, and with locally \ $L^{2}$ \ growth near \ $0$.\\
   It is easy to see that \ $\Xi_{\Lambda,S}^{(k)}$ \ is also a free  $\mathcal{O}_S[[a]]-$sheaf with the same generators. Both structures of \ $\mathcal{O}_S[[b]]$ \ and \ $\mathcal{O}_S[[a]]-$modules define a structure of \ $\mathcal{A}_S-$module on \ $\Xi_{\Lambda,S}^{(k)}$.\\
   More generally, for \ $V$ \ a finite dimensional complex vector space, we shall consider the \ $\mathcal{A}_S-$module \ $\Xi_{\Lambda,S}^{(k)}\otimes_{\C} V$ \ where \ $a$ \ and \ $b$ \ act as the identity on \ $V$. This is again a  free finite type module on \ $\mathcal{O}_S[[a]]$ \ and \ $\mathcal{O}_S[[b]]$. $\hfill \square$\\
   
   Recall that in the absolute case (i.e. \ $S = \{pt\}$) a geometric (a,b)-module is naturally embedded in an asymptotic (a,b)-module of the type  \ $\Xi^{(N)}_{\Lambda}\otimes_{\C} V$, where \ $\Lambda$ \ is a finite subset of \ $ \mathbb{Q} \cap ]0,1]$ \  (so \ $exp(-2i\pi.\Lambda)$ \ contains the eigenvalues of the monodromy), $N$ \ is an integer (a bound for the nilpotency order of the monodromy) and \ $V$ \ is a $k-$dimensional \ $\C-$vector space with \ $k : = rk(E)$\ which is defined as \ $Hom_{\A}(E, \Xi)$, where \ $\Xi : = \oplus_{\lambda \in \mathbb{Q} \cap ]0,1]} \ \Xi_{\lambda}$ \ and \ $\Xi_{\lambda} : = \cup_{N \in \mathbb{N}} \ \Xi_{\lambda}^{(N)}$. \\
   Note that in the case of the Brieskorn module associated to an isolated singularity of a holomorphic germ  \ $f$ \ in \ $\C^{n+1}$, we have \ $k = \mu$ \ the Milnor number and \ $V$ \ is the homology with complex coefficients of the Milnor fiber of \ $f$. But for any geometric (a,b)-module \ $E$ \ the vector space \ $V = Hom_{\A}(E, \Xi)$ \ is endowed with a quasi-unipotent automorphism \ $T$ \ induced by the monodromy of \ $\Xi$ \ (which is given by \ $Log s \mapsto Log s +  2i\pi$). This monodromic vector space depends only on the saturation  of \ $E$ \ by \ $b^{-1}.a$. See [B.95].

 \parag{Notation} When we consider a sheaf  $\mathcal{F}$ \ on \ $\mathcal{Y}$,  which is a complex space over \ $S\times \{0\}$ \   via the map \ $f$, we shall write that \ $\mathcal{F}$ \ is a \ $\mathcal{O}_{S}, \mathcal{O}_{S}[[a]], \A^{0}, \A-$module for short, when it is respectively  a module over the sheaf of algebras  on \ $\mathcal{Y}$
  $$f^{-1}(\mathcal{O}_{S}), f^{-1}(\mathcal{O}_{S}[[a]]), f^{-1}(\A^{0}), f^{-1}(\A).$$

\parag{The complex \ $(K_{/S}^{\bullet}, d_{/S}^{\bullet})$} In our situation \ $(@)$, under the hypothesis \ $(H0)$, we define the complex \ $(\hat{\Omega}_{/S}^{\bullet}, d_{/S}^{\bullet})$ \ which is the formal completion in \ $f$ \ of the holomorphic \ $S-$relative de Rham complex. Recall that, by definition, the sheaf \ $\hat{\Omega}_{/S}^{\bullet}$ \ is defined by the presheaf on \ $\mathcal{Y}$, given by the projective limit 
$$ U \longrightarrow  \lim\limits_{\longrightarrow N} \big[ \Omega_{/S}^{\bullet}(U)\big/f^N.\Omega_{/S}^{\bullet}(U)\big] .$$
The action of \ $a$ \ on \ $\hat{\Omega}_{/S}^{\bullet}$ \ is given by multiplication by \ $f$. It defines a structure of \ $\mathcal{O}_{S}[[a]]-$module on this sheaf. The degree \ $1$ \  differential
\ $d_{/S}^{\bullet}$ \ is induced by the \ $S-$relative de Rham differential. \\
Our smoothness assumption for the map \ $\pi: \mathcal{X} \to S$ \ implies that this complex is a resolution of \ $\pi^{-1}(\mathcal{O}_S)_{\vert \mathcal{Y}}$ \ with \ $\mathcal{O}_S-$linear differentials.\\
Then we define the sub-complex of \ $\mathcal{O}_{S}[[a]]-$modules \ $(K_{/S}^{\bullet}, d_{/S}^{\bullet})$ \ where the subsheaf \ $K^p_{/S}$ \ is the kernel of the map\footnote{we shall note \ $K^p_{f/S}$ \ when we want to specify the function \ $f$ \ we consider.}
$$ d_{/S}f \wedge : \hat{\Omega}_{/S}^p \to \hat{\Omega}_{/S}^{p+1} $$
 given by the left  exterior multiplication by \ $d_{/S}f$.  We shall denote \ $E^p$ \ the \ $p-$th cohomology sheaf of this complex.  Now the \ $\mathcal{O}_S[[a]]-$linearity of the differential \ $d_{/S}$ \ on this sub-complex  shows that these cohomology sheaves  are \ $\mathcal{O}_S[[a]]-$modules in an obvious way. We shall  explain now that, via the inverse of  the Gauss-Manin connection, there exists a natural action of  \ $b$ \  on these cohomology sheaves in such a way that they become \ $\mathcal{A}_S^0-$modules.
 
 \parag{Action of \ $b$ \ on the cohomology sheaves \ $E^{\bullet}$}
 
 Consider now the exact sequence of complexes of \ $\mathcal{O}_S[[a]]-$modules on \ $\mathcal{Y}$ 
  \begin{equation*} 0 \to (K_{/S}^{\bullet}, d_{/S}^{\bullet}) \to (\hat{\Omega}_{/S}^{\bullet},  d_{/S}^{\bullet}) \overset{d_{/S}f\wedge}{\longrightarrow}  (I_{/S}^{\bullet}, d_{/S}^{\bullet})[+1]  \to 0. \tag{1}
 \end{equation*}
 where \ $I_{/S}^p : = d_{/S}f\wedge \hat{\Omega}^{p-1}_{/S} $.
  Note that \ $K_{/S}^0$ \ and \ $I_{/S}^0$ \ are zero and that \ $K_{/S}^{1} = I_{/S}^1 = \hat{\mathcal{O}}_{\mathcal{X}}.d_{/S}f $, thanks to \ $(H0)$ : an holomorphic function on \ $\mathcal{X}\setminus \mathcal{Z}$ \ extends holomorphically to \ $\mathcal{X}$. \\
 The natural inclusions \ $ i : I_{/S}^p \hookrightarrow K_{/S}^p$ \ for all \ $p \geq 0$ \ are compatible with the differential \ $d_{/S}$. This leads to an exact sequence of complexes of \ $\mathcal{O}_S-$modules
 \begin{equation*}
 0 \to (I_{/S}^{\bullet}, d^{\bullet}_{/S}) \to (K_{/S}^{\bullet}, d^{\bullet}_{/S}) \to ([K_{/S}/\ I_{/S}]^{\bullet}, d_{/S}^{\bullet}) \to 0 .\tag{2}
 \end{equation*}
 where \ $[K_{/S}/\ I_{/S}]^{\bullet}$ \ is a coherent \ $\mathcal{O}_{\mathcal{X}}-$module with support contained in \ $\mathcal{Z}$. So  there exists, thanks to the Nullstellensatz and the hypothesis \ $(H0)$, locally  on \ $\mathcal{Z}$, an integer \ $N$ \ such that \ $f^N.K_{/S}^p \subset I_{/S}^p$ \ for any \ $p \geq 0$; so  we have 
 \begin{equation*}
f^{N}.[K_{/S}/\ I_{/S}]^p = 0  \quad \forall p \geq 0  \tag{3}
 \end{equation*}
When we shall assume that \ $\mathcal{Z}$ \ is \ $S-$proper, the direct image theorem of H. Grauert will give the \ $\mathcal{O}_S-$coherence of the direct images \ $R^q\pi_*( [K_{/S}\big/I_{/S}]^p )$. This will be a key point for our finiteness theorem (see section 2).\\
  We have a natural inclusion \ $f^{-1}(\hat{\Omega}^1_{S\times D/S}) \subset K_{/S}^1\cap Ker\, d_{/S}$, and this gives a sub-complex  concentrated in degree \ $1$ with zero differential of \ $(K_{/S}^{\bullet}, d_{/S}^{\bullet})$. As in [B.II], we shall consider also the quotient complex \ $(\tilde{K}_{/S}^{\bullet}, d_{/S}^{\bullet})$. So we have the exact sequence
 \begin{equation*}
  0 \to f^{-1}(\hat{\Omega}_{S\times D/S}^1) \to (K_{/S}^{\bullet}, d_{/S}^{\bullet}) \to (\tilde{K}_{/S}^{\bullet}, d_{/S}^{\bullet}) \to 0 . \tag{4}
  \end{equation*}
  This exact sequence corresponds to the (a,b)-module version of the exact sequence involving the nearby cycles and the vanishing cycles sheaves in the usual context (see for instance [B.II]  2, or more simply look for a fixed \ $s \in S$).\\
  Note that  with our hypothesis \ $(H0)$ \ the sheaf  \   $K_{/S}^1 \cap Ker\,d_{/S} = I^1_{/S} \cap Ker\, d_{/S}$ \ is isomorphic\footnote{If \ $g.d_{/S}f$ \ is \ $d_{/S}-$closed, it is a section of  \ $f^{-1}(\mathcal{O}_{S}[[a]])$ \ outside \ $\mathcal{Z}$, and it extends to \ $\mathcal{Y}$ \ thanks to \ $(H0)$.} to the  sheaf \ $f^{-1}(\mathcal{O}_{S}[[a]]).d_{/S}f$ \ on \ $\mathcal{Y}$ \  and the cohomology sheaves of  the complex \ $(\tilde{K}_{/S}^{\bullet}, d_{/S}^{\bullet})$ \ are supported in \ $\mathcal{Z}$. \\
  
  We want to describe now the natural actions of \ $a$ \ and \ $b$ \ on the cohomology sheaves of the complexes \ $(K^{\bullet}_{/S}, d_{/S}^{\bullet}), (I^{\bullet}_{/S}, d_{/S}^{\bullet}), ([K_{/S}/I_{/S}]^{\bullet}, d_{/S}^{\bullet}), f^*(\hat{\Omega}^1_{S\times D/S}), (\tilde{K}_{/S}^{\bullet}, d_{/S}^{\bullet})$.\\
  As in all cases we have complexes of \ $\mathcal{O}_S[[a]]-$modules, with  \ $\mathcal{O}_S[[a]]-$linear differentials,  the action of \ $a$ \ on the cohomology sheaves of these complexes is clear.\\
  
  The next lemma will give the action of \ $b$.
  
  \begin{lemma}\label{Tilde b}
  The exact sequence of complexes (1) induces  for any \ $p \geq 2$ \  a \ $\mathcal{O}_S-$linear isomorphism
   $$ \partial^p : \mathcal{H}^p(I_{/S}^{\bullet}, d_{/S}^{\bullet}) \to \mathcal{H}^p(K_{/S}^{\bullet}, d_{/S}^{\bullet})$$
   and a short exact sequence of \ $\mathcal{O}_S-$modules on \ $\mathcal{Y}$
   \begin{equation*}
    0 \to \mathcal{O}_{S}.d_{/S}f  \to \mathcal{H}^1(I_{/S}^{\bullet}, d_{/S}^{\bullet}) \overset{\partial^1}{\to} \mathcal{H}^1(K_{/S}^{\bullet}, d_{/S}^{\bullet}) \to 0 \tag{@@}
    \end{equation*}
    There is a canonical \ $\mathcal{O}_S-$linear splitting of \ $\partial^1$ 
    $$ \tilde{b}^1 :  \mathcal{H}^1(K_{/S}^{\bullet}, d_{/S}^{\bullet}) \to  \mathcal{H}^1(I_{/S}^{\bullet}, d_{/S}^{\bullet}) $$
    For \ $p \geq 2$ \  define \ $\tilde{b}^p : = (\partial^p)^{-1}$. Then we have for each \ $p \geq 0$ \  the formula
    \begin{equation*}
      a.\tilde{b} = \tilde{b}\circ\mathcal{H}^p( i)\circ\tilde{b} + \tilde{b}.a                             \tag{F}
     \end{equation*}
     where \ $i : (I_{/S}^{\bullet}, d_{/S}^{\bullet}) \to (K_{/S}^{\bullet}, d_{/S}^{\bullet})$ \ is the obvious map of \ $\mathcal{O}_S[[a]]-$complexes.
     \end{lemma}
 
 \parag{proof} The fact that for \ $p \geq 2$ \ the connector \ $\partial^p$ \ is an isomorphism is clear, as the \ $S-$relative de Rham complex is a \ $f^{-1}(\mathcal{O}_S)-$linear resolution of \ $f^{-1}(\mathcal{O}_S)$. \\
  Now we shall construct the splitting \ $\tilde{b}^1$. Consider \ $x \in K_{/S}^1$ \ such that \ $d_{/S}x = 0$. Thanks to the relative de Rham lemma we may write \ $x = d_{/S}\xi$ \ with \ $\xi \in \mathcal{O}_{\mathcal{X}}$. But the assumption that \ $x$ \ is \ $d_{/S}f\wedge-$closed implies that \ $\xi $ \ is locally constant along the \ $S-$fibers of\ $\mathcal{Y} \setminus \mathcal{Z}$ \ because  near such a point, we may find  \ $S-$relative coordinates \ $(z_0, \dots, z_n)$ \  on \ $\mathcal{X}$ \ such that \ $f(s,z) = z_0$. Then the relation \ $d_{/S}\xi\wedge d_{/S}f = 0$ \ implies that \ $\xi(s,z)$ \ is a function of \ $(s,z_0)$ \ only. This proves our assertion, and the fact that \ $\mathcal{Y}$ \ is locally connected allows to choose \ $\xi$ \ vanishing (set theoretically) on \ $\mathcal{Y}$, because  \ $\mathcal{Y} \setminus \mathcal{Z}$ \ is dense in \ $\mathcal{Y}$. Moreover this choice is unique because in the previous  local computation near a  point in \ $\mathcal{Y} \setminus \mathcal{Z}$, if \ $x = g(s,z_0).dz_0$ \ we have to choose \ $\xi = G(s,z_0)$ \ with \ $\frac{\partial G}{\partial z_0}(s, z_0) = g(s,z_0)$ \ and \ $G(s,0) = 0$, which has an unique solution. So the continuity of \ $\xi$ \ forces vanishing on all \ $\mathcal{Y}$.\\
 Now we define\footnote{Note that \ $d_{/S}\xi \wedge d_{/S}f = 0$ \ implies that \ $d_{/S}f\wedge \xi$ \ is \ $d_{/S}-$closed.} \ $\tilde{b}^1[d_{/S}\xi] = [d_{/S}f\wedge \xi]$ \ with  this unique choice of \ $\xi$. It is clearly \ $\mathcal{O}_S-$linear and 
 $$\partial^1\circ \tilde{b}^1[d_{/S}\xi] = \partial^1 [d_{/S}f\wedge \xi] = [d_{/S}\xi] $$
 as \ $\partial^1 = d_{/S}\circ \big(d_{/S}f\wedge\big)^{-1}$.\\
   To finish the proof we shall prove the formula \ $(F)$ \ for each \ $p \geq 0$. Writing  \ $x = d_{/S}\xi$, we get
   $$ a[x] + \mathcal{H}^p(i)[\tilde{b}[x]] =[ f.d_{/S}\xi + d_{/S}f \wedge \xi] = [d_{/S}(f.\xi)] $$
   and so
   $$ \tilde{b}( a[x] + \mathcal{H}^p(i)[\tilde{b}[x]]) = [ d_{/S}f \wedge f.\xi ] = a.\tilde{b}([x]) $$
   which concludes the proof. $\hfill \blacksquare$
   
   \begin{defn}\label{b en cohomologie}
   We shall define
   \begin{enumerate}
   \item  \ $b : \mathcal{H}^p(K^{\bullet}_{/S}, d_{/S}^{\bullet}) \to \mathcal{H}^p(K^{\bullet}_{/S}, d_{/S}^{\bullet})$ \ as \ $ b : = \mathcal{H}^p(i)\circ \tilde{b}$.
   \item \ $ b : \mathcal{H}^p(I^{\bullet}_{/S}, d_{/S}^{\bullet}) \to \mathcal{H}^p(I^{\bullet}_{/S}, d_{/S}^{\bullet})$ \ as \ $ b : = \tilde{b}\circ \mathcal{H}^p(i) $.
   \item \ $ b = 0$ \ on the cohomology of the quotient complex \ $ ([K_{/S}/I_{/S}]^{\bullet}, d_{/S}^{\bullet})$.
   \item \ $ b : = \int_0^z $ \ on \ $ f^{-1}(\hat{\Omega}^1_{S\times D/S}) \simeq \mathcal{O}_S[[z]].dz $.
   \item  and \ $b$ \ is induced on \ $\mathcal{H}^p(\tilde{K}_{f/S}, d_{/S})$ \ by \ $b$ \ on \ $\mathcal{H}^p(K^{\bullet}_{/S}, d_{/S}^{\bullet})$.
     \end{enumerate}
   \end{defn}
   
 \parag{Remarks} \begin{enumerate}
 \item   For \ $p \not= 1$ \ the cohomology sheaves  \ $\mathcal{H}^p(K^{\bullet}_{/S}, d_{/S}^{\bullet})$ \ and  \ $\mathcal{H}^p(\tilde{K}_{f/S}, d_{/S})$ \ are canonically isomorphic as \ $\A^{0}-$modules  and that for \ $p = 1$ \   the maps in the exact sequence
   $$ 0 \to  f^*(\hat{\Omega}^1_{S\times D/S}) \to \mathcal{H}^1(K^{\bullet}_{/S}, d_{/S}^{\bullet}) \to \mathcal{H}^1(\tilde{K}_{f/S}, d_{/S}) \to 0 $$
   are  \ $\A^{0}-$linear.
  \item  Note that under the hypothesis \ $(H0)$ \ we have the vanishing of \ $\mathcal{H}^1(\tilde{K}_{f/S}, d_{/S})$.\\
  \end{enumerate}  
   
   An easy consequence of the previous lemma is that we have \ $a.b - b.a = b^2$ \ for each cohomology sheaf of any of these complexes of \ $\mathcal{O}_S[[a]]-$modules.

 \parag{the complexes \ $(\mathcal{K}_{/S}^{0,\bullet}, D_{/S}^{\bullet})$ \ and $(\mathcal{K}_{/S}^{\bullet}, D_{/S}^{\bullet})$} Define first the subsheaves \\
  $\mathcal{K}_{/S}^{\bullet} \subset  \hat{\Omega}_{/S}^{\bullet}[[b]]  $ \  by the condition 
  
 $$ \Omega : =  \sum_{j\geq 0}^{+\infty} \ b^j.\omega_j \in \mathcal{K}_{/S}^{\bullet} \quad {\rm when} \quad   d_{/S}f\wedge \omega_0  = 0 .$$
 
  Then define \ $ \mathcal{K}_{/S}^{0,\bullet} = \hat{\Omega}_{/S}^{\bullet}[b]  \cap \mathcal{K}_{/S}^{\bullet} $. The differential \ $D^{\bullet}_{/S} $ \ is given by 
  $$ D^{\bullet}_{/S}(\sum_{j=0}^{+\infty} \ b^j.\omega_j) = \sum_{j=0}^{+\infty} \ b^j.(d_{/S}\omega_j - d_{/S}f\wedge \omega_{j+1}) .$$
  It sends \ $\mathcal{K}_{/S}^p$ \ to \ $\mathcal{K}_{/S}^{p+1}$ \ and also \  \ $\mathcal{K}_{/S}^{0,p}$ \ to \ $\mathcal{K}_{/S}^{0,p+1}$. Moreover it satisfies \ $D^{p+1}_{/S}\circ D^p_{/S} = 0 $ \ for each \ $p \geq 0$.\\
   Now define the action of \ $a$ \ and \ $b$ \ as follows :
           \begin{align*}
     & a.\sum_{j=0}^{+\infty} b^j.\omega_j  = \sum_{j=0}^{+\infty} b^j.(f.\omega_j + (j-1).\omega_{j-1}) \quad {\rm with \ the \ convention} \quad \omega_{-1} = 0 \\
     & b.\sum_{j=0}^{+\infty} b^j.\omega_j  = \sum_{j=1}^{+\infty} b^j.\omega_{j -1}
     \end{align*}
     It is clear that \ $\mathcal{K}_{/S}^{0,\bullet}$ \ is a \ $\mathcal{O}_{/S}[a]-$module, but the fact that 
   \ $\mathcal{K}^{\bullet}_{/S}$ \ is a  \ $\mathcal{O}_S[[a]]-$module is a consequence of the following lemma.
   
   \begin{lemma}\label{action a}
   Let \ $\Omega \in \mathcal{K}_{/S}^{\bullet}$. Then for each \ $N \in \mathbb{N}$ \ and \ each \ $(h,j) \in \mathbb{N}^2, h \leq j$,  there exists a polynomial \ $T^N_{j,h}(a)$ \ of degree \ $\leq N$ \ and of valuation at \ $0$ \ at least equal to \ $ N-h$ \ in \ $\C[a]$,  such that the \ $j-$th component of \ $a^N.\Omega$ \ is given by
   \begin{equation*}
    (a^N.\Omega)_j = \sum_{h=0}^j \ T^N_{j,h}(a).\omega_{j-h} . \tag{*}
    \end{equation*}
   \end{lemma}
   
   \parag{proof} We shall prove the relation \ $(*)$ \ and the estimations on the degree and on the valuation at \ $0$ \  of \ $T^N_{j,h}$ \ by induction on \ $N$. We have
   \begin{align*}
   & (a^{N+1}.\Omega)_j = f.(a^N.\Omega)_j + (j-1).(a^N.\Omega)_{j-1} \\
   & \qquad  = \sum_{h=0}^j \big(a.T^N_{j,h}(a) + (j-1).T^N_{j-1,h-1}(a)\big).\omega_{j-h} 
   \end{align*}
   with the convention \ $T^N_{j,h} = 0$ \ if \ $j < 0$ \ or \ $h < 0$. This gives the induction relation
   $$ T^{N+1}_{j,h}(a) = a.T^N_{j,h}(a) + (j-1).T^N_{j-1,h-1}(a) .$$
   The assertions on the degree and on the valuation are immediate. $\hfill \blacksquare$\\
   
   Now the verification of the  following facts are easy
   \begin{enumerate}[i)]
   \item \ $\mathcal{K}^{\bullet}_{/S}$ \ is a \ $\mathcal{A}_S-$module.
   \item The differential \ $D^{\bullet}_{/S} $ \ is \ $\mathcal{A}_S-$linear.
   \end{enumerate}
   
   \parag{Notation} We define \ $\mathcal{E}^{p}$ \ as the \ $p-$th cohomology sheaf of the complex \ $(\mathcal{K}_{/S}^{\bullet}, D_{/S}^{\bullet})$. They are sheaves of \ $\mathcal{A}_S-$modules on \ $\mathcal{Y}$.
   
   \begin{lemma}\label{image de b}
   Let \ $X : = \sum_{j=0}^{\infty} \ b^j.x_j$ \ be a \ $D_{/S}-$closed element in \ $\mathcal{K}_{/S}^p$. Then \ $[X] $ \ is in \ $b.\mathcal{E}^p$ \ if and only if \ $x_0$ \ is in \ $I_{/S}^p + d_{/S}K^{p-1}_{/S}$.
   \end{lemma}
   
   \parag{proof} A relation  \ $X = b.Y + D_{/S}U $ \ with \ $Y \in \mathcal{K}_{/S}^p$, (and also $D_{/S}-$closed), and \ $U \in \mathcal{K}_{/S}^{p-1}$ \ implies \ $x_0 = d_{/S}(u_0) - d_{/S}f\wedge u_1$, so our assertion is clear. \\
   Conversely, if \ $x_0 = d_{/S}f\wedge w + d_{/S}(v)$ \ for some \ $w \in \hat{\Omega}^{p-1}_{/S}$ \ and some \ $v$ \ in \ $ K^{p-1}_{/S}$ \  we have \ $D_{/S}(v - b.w) =  d_{/S}f\wedge w  + d_{/S}(v) - b.d_{/S}(w)$ \ and then we have  \ $ [X] = [X -  D_{/S}(v- b.w)] = b.[Y]$ \ where \ $Y : = x_1 + d_{/S}(w) + \sum_{j=2}^{\infty} \ b^{j-1}.x_j$. And \ $ D_{/S}(b.Y) = 0$ \ implies \ $d_{/S}f\wedge (x_1 +  d_{/S}(w)) = 0$, so \ $Y $ \ is in \ $\mathcal{K}^{p-1}_{/S}$. $ \hfill \blacksquare$\\
   
   \begin{cor}\label{inj. coh.}
   The natural map \ $\mathcal{K}_{/S}^{\bullet} \to K_{/S}^{\bullet}$ \ defined by \ $X : = \sum_{j=0}^{\infty} \ b^j.x_j \mapsto x_{0}$ \ induces an injective \ $\mathcal{O}_S[[a]]-$linear map \ $\mathcal{E}^p\big/b.\mathcal{E}^p \to K^p_{/S}\big/(I^p_{/S} + d_{/S}K_{/S}^{p-1})$. 
   \end{cor}
   
   \parag{proof} First remark that if \ $X$ \ is in \ $\mathcal{K}_{/S}^{p}$ \ and \ $U$ \ in \ $\mathcal{K}_{/S}^{p-1} $ \ then the \ $0-$component of \ $X + D_{/S}U$ \ is equal to \ $x_{0} + d_{/S}u_{0} -d_{/S}f\wedge u_{1}$ \ so gives the same class than \ $x_{0}$ \ in the quotient \ $K^p_{/S}\big/(I^p_{/S} + d_{/S}K_{/S}^{p-1})$. Also the \ $0-$th component of \ $a.X$ \ is \ $a.x_{0}$ \ for the natural action of \ $\mathcal{O}_{/S}[[a]]$ \ on \ $K^p_{/S}$. As \ $I^p_{/S} + d_{/S}K_{/S}^{p-1}$ \ is a \ $\mathcal{O}_{/S}[[a]]-$submodule of \ $K^p_{/S}$ \ the map \ $\mathcal{E}^{p} \to K^p_{/S}\big/(I^p_{/S} + d_{/S}K_{/S}^{p-1})$ \ is well defined and \ $\mathcal{O}_{/S}[[a]]-$linear. It is clear that this map vanishes on \ $b.\mathcal{E}^p$ \ and the previous lemma gives the injectivity of the induced map \ $\mathcal{E}^p\big/b.\mathcal{E}^p \to K^p_{/S}\big/(I^p_{/S} + d_{/S}K_{/S}^{p-1})$.$\hfill \blacksquare$\\

An  easy  consequence of this corollary is that,  using the hypothesis \ $(H0)$ \ and the Nullstellensatz which gives locally on \ $\mathcal{Y}$ \ an integer \ $N$ \ such that \ $f^{N}.K_{/S}^{\bullet} \subset I_{/S}^{\bullet}$ , we may find  locally on \ $\mathcal{Y}$ \ an integer \ $N$ \ such that \ $a^N.\mathcal{E}^p \subset b.\mathcal{E}^p \quad \forall p \geq 0$.\\
   
   \begin{prop}[quasi-isomorphism \ $u^0$.]\label{quasi-iso}
The natural map given by the inclusion  \ $u^0(x) = x$ \  is a quasi-isomorphism \ $u^0 : (K_{/S}^{\bullet}, d_{/S}^{\bullet}) \to (\mathcal{K}^{0,\bullet}_{/S}, D_{/S}^{\bullet})$ \ of complexes of \ $\mathcal{O}_{/S}[[a]]-$modules. It preserves also the action of \ $b$ \ on the cohomology sheaves, and so induces isomorphisms of \ $\mathcal{A}^0_S-$modules on the cohomology sheaves.
\end{prop}

\parag{proof} The identities \ $D_{/S}\circ u^0 = u^0\circ d_{/S}$ \ and \ $a\circ u^0 = u^0\circ a$ \ are easily verified. For the action of \ $b$ \ on cohomology, take \ $x \in K^p_{/S}$ \ such that \ $d_{/S}(x) = 0$. Then writing \ $x = d_{/S}(\xi)$ \ for some \ $\xi \in \hat{\Omega}^{p-1}_{/S}$ \ (recall that  the complex \ $( \hat{\Omega}^{\bullet}_{/S}, d_{/S}^{\bullet})$ \ is a resolution of \ $f^{-1}(\mathcal{O}_{S})$ \ on \ $\mathcal{Y}$), we have \ $b.[x] = [d_{/S}f\wedge \xi]$. Now we have \ $b.\xi \in \mathcal{K}^{0,p-1}_{/S}$ \ and \ $D_{/S}(b.\xi) = -d_{/S}f\wedge \xi + b.d_{/S}(\xi) $. So we have \ $ D_{/S}(b.\xi) = -d_{/S}f\wedge \xi + x.b$ \ in \ $\mathcal{K}^{0,p}_{/S}$. This gives \ $b\circ u^0 = u^0\circ b$ \ in the cohomology sheaf \ $\mathcal{H}^p(\mathcal{K}^{0,\bullet}_{/S}, D_{/S})$.\\
We shall prove now the injectivity of the map \ $\mathcal{H}^p(u^0)$. Let \ $x$ \ be in \ $K^p_{/S}$ \ such that \ $d_{/S}(x) = 0$ \ and assume that \ $u^0(x) = D_{/S}Y$ \ where \ $Y : = \sum_{j = 0}^N \ b^j.y_j$ \ with \ $y_0 \wedge d_{/S}f = 0 $. Then we have
\begin{equation*}
 x = d_{/S}(y_0) - d_{/S}f \wedge y_1 \qquad {\rm and} \quad d_{/S}(y_j) = d_{/S}f\wedge y_{j+1} \quad \forall j \geq 1 \tag{@}
\end{equation*}
For \ $j=N$ \ we have \ $d_{/S}(y_N) = 0 $, so we may write (locally) \ $y_N = d_{/S}(\eta_N)$ \ for some \ $\eta_N \in \hat{\Omega}^{p-1}_{/S}$. Now we have \ $d_{/S}(y_{N-1}) = d_{/S}f\wedge y_N = -d_{/S}(d_{/S}f\wedge \eta_N)$;  so we obtain \ $d_{/S}(y_{N-1} + d_{/S}f\wedge \eta_N) = 0 $. Then going on this way, we construct \ $\eta_j \in  \hat{\Omega}^{p-1}_{/S}$ \ such that
\begin{equation*}
 y_j =  d_{/S}(\eta_j) - d_{/S}f\wedge \eta_{j+1} \qquad \forall j \in [1,N]  \tag{@@}
 \end{equation*}
So we get \ $x = d_{/S}(y_0) + d_{/S}(d_{/S}f\wedge\eta_1)$ \ and this shows that \ $[x] = 0$ \ in the cohomology of \ $(K_{/S}^{\bullet}, d_{/S}^{\bullet})$.\\
To prove the surjectivity on cohomology, take \ $Y \in \mathcal{K}^{0,p}_{/S}$ \ such that \ $D_{/S}Y = 0$. Write \ $Y = \sum_{j=0}^N \ b^j.y_j$; then the second relation \ $(@)$ \ is satisfied for each \ $j \geq 0$ \ and arguing as above we may find \ $\eta_j \in \hat{\Omega}^{p-1}_{/S}$ \ such that \ $(@@)$ \ holds for each \ $j \geq 1$ \ (with \ $\eta_j = 0$ \ for \ $j \geq N+1$). Then let \ $H : = \sum_{j=1}^N \ b^j.\eta_j$. We obtain
$$ Y - D_{/S}H = y_0 + d_{/S}f\wedge \eta_1  .$$
But  \ $y_0 + d_{/S}f\wedge \eta_1$ \ is in \ $K^p_{/S} \cap Ker\, d_{/S}$ \ as \ $d_{/S}(y_0) = d_{/S}f\wedge d_{/S}(\eta_1)$. So we have \ $\mathcal{H}^p(u^0)[y_0 + d_{/S}f\wedge \eta_1] = [Y]$. $\hfill \blacksquare$

\parag{Notation} We shall identify via \ $\mathcal{H}(u^0)$ \ the cohomology sheaf \ $E^p$ \ of the complex  \ $(K_{/S}^{\bullet}, d_{/S}^{\bullet})$ \ to the cohomology sheaf of the complex of \ $\mathcal{A}^0_S-$modules \ $ (\mathcal{K}^{0,\bullet}_{/S}, D_{/S})$. 

\begin{prop}[the map \ $u^1$.]\label{u1}
The obvious map \ $u^1 : (\mathcal{K}^{0,\bullet}_{/S}, D_{/S}^{\bullet}) \to  (\mathcal{K}^{\bullet}_{/S}, D_{f/S}^{\bullet})$ \ induces  a map \ $\mathcal{H}^p(u^1)$ \ with kernel contained in \ $\cap_{m \geq 0} \ b^m.E^p$.\\
 If we have  \ $\cap_{m \geq 0} \  b^m(E^p) = \{0\} $ \ then the  image of \ $\mathcal{H}^p(u^1)$  \  is dense  for the \ $b-$adic topology of \ $\mathcal{E}^p$,  the \ $p-$th cohomology sheaf of the complex \ $ (\mathcal{K}^{\bullet}_{/S}, D_{/S}^{\bullet})$. 
\end{prop}

Note that we dont know that the \ $b-$adic topology of  \ $\mathcal{E}^p$ \ is separated under the hypothesis of this proposition (which is only \ $(H0)$ \ as in all this section 1.2).

\parag{Proof} We begin by looking at  the kernel of \ $\mathcal{H}^p(u^1)$. Using the quasi-isomorphism \ $u^0$ \ we may consider only some \ $x_0 \in K^p_{/S} \cap Ker\, d_{/S}$ \ such that there exists an  \ $Y = \sum_{j=0}^{\infty} b^j.y_j $ \ in \ $ \mathcal{K}^p_{/S}$ \ such that \ $D_{/S}Y = 0$. Then we have
$$ x_0 = d_{/S}(y_0) - d_{/S}f\wedge y_1 \quad {\rm and} \quad d_{/S}(y_j) = d_{/S}f \wedge y_{j+1} \quad  \forall j\geq 1 .$$
As \ $[d_{/S}(y_j)] = b.[d_{/S}(y_{j+1})] $ \ for each \ $j \geq 1$ \ we obtain that \ $[x_0] = [x_0 - d_{/S}(y_0)]$ \ is in \ $b^m(E^{p})$ \ for each \ $m \geq 0$. This proves the first point.\\
To see that the image is dense, let \ $X \in \mathcal{K}^p_{/S}$ \ such that \ $D_{/S}X= 0$ \ and fix an integer \ $N$. Let \ $X = \sum_{j=0}^{\infty} \ b^j.x_j$ \ and define \ $X_N : = \sum_{j=0}^N \ b^j.x_j$. As we have  \ $d_{/S}(x_j) = d_{/S}f\wedge x_{j+1} \quad \forall j\geq 0$, we have \ $D_{/S}(X_N) = d_{/S}(x_N).b^N $. But the relation \ $d_{/S}(x_j) = d_{/S}f\wedge x_{j+1}$ \ for all \ $j \geq 0$ \ implies that \ $[d_{/S}(x_N)] \in \cap_{m \geq 0}\ b^m.E^{p}$;  so, thanks to our hypothesis \ $\cap_{m \geq 0} \  b^m(E^p) = \{0\} $, we may find \ $z_N$ \ such that \ $d_{/S}(z_N) = d_{/S}(x_N) $ and \ $d_{/S}f \wedge z_N = 0$. We have \ $D_{/S}(b^N.z_N) = b^N.d_{/S}(x_N)$ \ and so \ $Y_N : = X_N - b^N.z_N$ \ is \ $D_{/S}-$closed and satisfies \ $X - Y_N \in b^N.\mathcal{K}_{/S}^p$. As \ $Y_N : = X_N - b^N.z_N$ \ is in \ $\mathcal{K}_{/S}^{0,p}$ \ we have approximated \ $X$ \ modulo \ $b^{N}.\mathcal{E}^{p}$ \ by an element in the image of \ $\mathcal{H}^p(u^1)$.$\hfill \blacksquare$\\

\section{The quasi-isomorphism.}

We shall begin by showing how we shall use our hypotheses \ $(H1)$ \ and \ $(H2)$.

\subsection{Use of \ $(H1)$.}

The hypothesis \ $(H1)$ \ will be  useful thanks to the following proposition.

 \begin{prop}\label{a-torsion}
Let \ $\mathcal{X} \overset{\pi}{\longrightarrow} S$ \ be a \ $S-$relative complex manifold and let \\
 $f : \mathcal{X} \to D $ \ be an holomorphic function satisfying the following assumptions :
\begin{enumerate}[i)]
\item \ $\mathcal{Z}  : = \{d_{/S}f = 0\}  \subset \mathcal{Y} : =  f^{-1}(0)$.
\item Locally on \ $\mathcal{X}$ \ the function \ $f$ \ admits a \ $S-$desingularization.
\end{enumerate}
Then there exists locally on \ $\mathcal{X}$ \ an integer  \ $N$ \ such that we have \ $a^N.A(E^{\bullet}) = 0$ \ where \ $E^{\bullet}$ \ are the cohomology sheaves of the complex  \ $(K^{\bullet}_{f/S}, d_{/S})$.
\end{prop} 

First recall the following lemma (for a proof see  for instance [B-S 04] prop. 2.2.)

\begin{lemma}\label{cas S-DCN}
We consider the situation \ $(@)$ \ where we assume that \ $\mathcal{X} : = S \times U$ \ where \ $U \subset \C^{n+1}$ \ is a open neighbourghood of the origin and where \ $\pi$ \ is the first projection. We assume that \ $f(s,z) = z^{\alpha}$ \ where \ $(z_0, \dots,z_n)$ \ are the coordinates on \ $\C^{n+1}$ \ and where \ $\alpha \in \mathbb{N}^{n+1}$. Then the \ $a-$torsion \ $A(E^{\bullet})$ \ of the cohomology sheaves of the complex \ $(K_{/S}^{\bullet}, d_{/S}^{\bullet})$ \ is \ $0$ \ and \ $b^{-1}.a$ \ acts bijectively on the cohomology sheaves \ $E^{\bullet}$.
\end{lemma}

\parag{proof of the proposition} As the problem is local on \ $\mathcal{X}$ \  we may assume the following facts :
\begin{enumerate}[i)]
\item \ $\mathcal{X} = S \times U$ \ where \ $S$ \ is Stein and where \ $U$ \ is a polydisc in \ $\C^{n+1}$.
\item We have a \ $S-$proper modification \ $\tau : \tilde{\mathcal{X}} \to \mathcal{X}$ \ such that \ $ \tilde{\mathcal{X}}$ \ is \ $S-$smooth and the zero set of \ $\tilde{f} : = f\circ \tau$ \ defines  on \ $\tilde{\mathcal{X}}$ \ a \ $S-$relative normal crossing divisor.
\end{enumerate}
We have, for each open set \ $\mathcal{X}'$ \ in \ $\mathcal{X}$ \ and for each integer \ $q$,  an "edge" map from the spectral sequence of the hypercohomology, 
\begin{equation*}
\varepsilon_{q} : H^q\big(\Gamma(\mathcal{X}', K^{\bullet}_{f/S}), d^{\bullet}_{/S}\big) \to \mathbb{H}^q\big(\tilde{\mathcal{X}}', (K^{\bullet}_{\tilde{f}/S}, d_{/S}^{\bullet})\big) \tag{@}
\end{equation*}
where \ $\tilde{\mathcal{X}}' = \tau^{-1}(\mathcal{X}')$. It is consequence of the equalities \ $\tau_*(K^{\bullet}_{\tilde{f}/S}) = K^{\bullet}_{f/S}$ \ which are compatible with the \ $S-$relative de Rham differential. These maps define  morphisms, denoted \ $e_{q}$,  of the corresponding   \ $\mathcal{O}_S[[a]]-$sheaves on \ $\mathcal{X}$. 

\parag{Claim} There exists locally on \ $\mathcal{X}$ \ an integer \ $N_1$ \ such the sheaf maps \ $e_{q}$ \ induce  isomorphisms between \ $a^{N_1}.\mathcal{H}^q\big(K^{\bullet}_{f/S}, d_{/S}^{\bullet}\big)$ \ onto \ $a^{N_1}. \mathcal{H}^q\big(\mathbb{R}\tau_*(K^{\bullet}_{\tilde{f}/S}, d_{/S})\big) $ \ where \ $\mathbb{R}\tau_*(K^{\bullet}_{\tilde{f}/S}, d_{/S}^{\bullet})$ \ is the direct image in the category of complexes of sheaves.

\parag{proof of the claim} First remark that for each \ $i \geq 1$ \ the sheaves \ $R^i\tau_*(K^{\bullet}_{\tilde{f}/S})$ \ are \ $a-$completions of  coherent sheaves  on \ $\mathcal{X}$ \ and supported in \ $\mathcal{Y} : = f^{-1}(0)$. So on any relatively compact Stein open set in \ $\mathcal{Y}$ \  we may find an integer \ $N_1$ \ such that \ $f^{N_1}.R^i\tau_*(K^{\bullet}_{\tilde{f}/S}) = 0 $. This means that, in the spectral sequence
$$ E_2^{p,q} : = \mathcal{H}^q\big(R^p\tau_*(K^{\bullet}_{\tilde{f}/S}), d^{\bullet}_{/S}\big) $$
which converges to \ $\mathcal{H}^{m}\big(\mathbb{R}\tau_*(K^{\bullet}_{\tilde{f}/S}, d_{/S}^{\bullet})\big)$,  the sub-spectral sequence \ $a^{N_1}.E_2^{p,q}$ \ which converges to \ $a^{N_1}.\mathcal{H}^{m}\big(\mathbb{R}\tau_*(K^{\bullet}_{\tilde{f}/S}, d_{/S}^{\bullet})\big)$, degenerates at \ $E_2$ \ and gives that each  \ $e_{q}$ \ induces the isomorphism of the claim. $\hfill \blacksquare$

\parag{end of the proof of \ref{a-torsion}} Now it is enough to prove that there exists locally on \ $\mathcal{X}$ \ an integer \ $N_2$ \ such that if \ $x \in  \mathcal{H}^q\big(K^{\bullet}_{f/S}, d_{/S}^{\bullet}\big)$ \ is in the \ $a-$torsion of this sheaf, then its image in \ $ \mathcal{H}^q\big(\mathbb{R}\tau_*(K^{\bullet}_{\tilde{f}/S}, d_{/S}^{\bullet})\big) $ \ is killed by \ $a^{N_2}$ \ in this sheaf, because this allows to conclude by taking \ $N : = N_1+N_2$.\\
This is a easy exercice on hypercohomology because at each step of the spectral sequence we know that \ $x$ \ is locally killed on \ $\tilde{\mathcal{X}}$ \ by \ $a$,  thanks to the lemma \ref{cas S-DCN},  and that the obstruction for the next step is in the sheaf \ $R^i\tau_*(K^{q-i}_{\tilde{f}/S})$, which is killed by \ $a^{N_{1}}$. So we may take \ $N_2 = (n+1).(N_1+1)$. $\hfill \blacksquare$

 \subsection{Use of \ $(H2)$.}
 
 First recall a simple lemma from [B.06] lemma 2.1.2
 
 \begin{lemma}\label{a donne b}
The commutation relation \ $a.b - b.a = b^2$ \ implies for each integer \ $N$ \  the formula
$$ N!.b^{2N} = \sum_{j=0}^N \ (-1)^j.\frac{N!}{(N-j)!.j!}.b^j.a^N.b^{N-j} .$$
\end{lemma}

It implies that an \ $\mathcal{A}^0_S-$module which is killed by \ $a^N$ \ is killed by \ $b^{2N}$. \\

We shall  use  the hypothesis \ $(H2)$ \ via the following proposition :

\begin{prop}\label{b-separation}
In the situation \ $(@)$ \ with the hypotheses \ $(H0), (H1)$ \ and \ $(H2)$ \ the cohomology sheaves \ $E^{\bullet}$ \  of the complex \ $(K_{/S}^{\bullet}, d_{/S}^{\bullet})$ \ are \ $b-$separated, that is to say satisfy \ $\cap_{m \geq 0} \ b^m.E^{\bullet} = \{0\}$.
\end{prop}

\parag{proof} We begin by proving that it is enough to establish the following  key properties of \ $E : = E^{\bullet}$ \  :
\begin{enumerate}[i)]
\item There exists locally on \ $\mathcal{Y}$ \ and integer \ $N$ \ such that \ $a^N.A(E) = \{0\}$.
\item \ $B(E) \subset A(E)$.
\item  \ $ \cap_{m \geq 0} \ b^m.E \subset A(E)$.
\end{enumerate}
Remark first that the inclusion \ $\cap_{m \geq 0} \ b^m.E^{\bullet} \subset B(E)$ \ implies \ $\cap_{m \geq 0} \ b^m.E^{\bullet} = \{0\}$. Moreover the inclusion ii) implies that \ $B(E) \subset \tilde{A}(E)$ \ because \ $\tilde{A}(E)$ \ is the biggest  \ $\mathcal{A}^0_S-$submodule in \ $A(E)$ \ (by definition), and \ $B(E)$ \ is \ $\mathcal{A}^0_S-$stable. But the inclusion  \ $B(E) \subset \tilde{A}(E)$ \  implies the equality \ $B(E) =  \tilde{A}(E)$ \  because on \ $\tilde{A}(E)\big/B(E)$ \ the map \  $b$ \ is injective and satisfies \ $ b^{2N} = 0$, thanks to the lemma \ref{a donne b} and the condition i). Now \ $\cap_{m \geq 0} \ b^m.E^{\bullet}$ \ is a \ $\mathcal{A}^0_S-$submodule, so condition iii) implies the inclusion \ $\cap_{m \geq 0} \ b^m.E^{\bullet} \subset \tilde{A}(E) = B(E)$ \ and our claim is proved.\\
The condition i) has been established in proposition \ref{a-torsion}. The conditions ii) and iii)  will be deduced from the following lemma.

\begin{lemma}\label{dev. asympt.}
In the situation \ $(@)$ \ with the hypothesis \ $(H0)$ \ and \ $(H2)$ \ let \ $y \in \mathcal{Y}$ \ and denote \ $F_{y}$ \ the Milnor fiber  of \ $f$ \ at \ $y$. Put \ $\pi(y) : = s_0$.  We have a natural \ $\big(\mathcal{O}_S[[a]]\big)_{s_0}-$linear map
$$ dev_{y} : E^{p+1}_{y} \to \big(\Xi_{\Lambda, S}^{(n)}\big)_{s_0} \otimes_{\C} H^p(F_{y},\C) $$
which commutes with the actions of \ $b$. The kernel of \ $dev_{y}$ \ is exactly the germ at \ $y$ \ of the  \ $a-$torsion \ $A(E^{p+1})$ \ of \ $E^{p+1}$.
\end{lemma}

\parag{proof} Let \ $[\omega] \in E^{p+1}_y$ \ where \ $\omega $ \ is a local section of \ $K_{/S}^{p+1}$ \ which is \ $d_{/S}-$closed, and, for \ $\gamma \in H_p(F_y, \C)$, let \ $(\gamma_{s,z}), (s,z) \in S\times D^*$ \ be the horizontal multivalued family of \ $p-$cycles in the fibers deduced form \ $\gamma$ ;  then define
$$ dev_y[\omega] : = \big[ \int_{\gamma_{s,z}} \ \omega\big/d_{/S}f\big] \in \big(\Xi_{\Lambda, S}^{(n)}\big)_{s_0} $$
where \ $\big(\Xi_{\Lambda, S}^{(n)}\big)_{s_0}$ \ is the germ at \ $s_0$ \ of multivalued asymptotic expansions (in the local coordinate \ $z$ \  in \ $D$, see the example in section 1.2) \  where \ $exp(2i\pi.\Lambda)$ \ is the set of eigenvalues of the local monodromy at \ $y$ \ and where \ $n = \dim_{\C} F_y$ \ bounds the degrees of the logarithms. This commutes with the action of \ $\mathcal{A}_{S,s_0}$. Then the fact that \ $[\omega]$ \ is in the kernel of \ $dev_y$ \ is clearly equivalent to the fact that \ $[\omega]$ \ induces the zero cohomology class on the generic fibers of the \ $S-$relative Milnor fibration of \ $f$ \ at \ $y$. So this exactly means that \ $[\omega]$ \ is in \ $A(E^{p+1})_y$. $\hfill \blacksquare$\\

\parag{End of the proof of the proposition \ref{b-separation}} Note that it is enough to prove properties ii) and iii) for \ $E^p_y$ \ for any \ $p \geq 0$ and any \ $y \in \mathcal{Y}$. But  they are immediate consequence of the fact that the \ $\mathcal{A}_{S,s_0}-$module \ $\big(\Xi_{\Lambda, S}^{(n)}\big)_{s_0} \otimes_{\C} H^p(F_{y},\C) $ \ has no \ $a-$torsion and is \ $b-$separated. $\hfill \blacksquare$\\

We have another important consequence of the lemma \ref{dev. asympt.}

\begin{cor}\label{completion 1}
In the situation \ $(@)$ \ with the hypothesis \ $(H0)$ \ and \ $(H2)$ \ the \ $(\mathcal{O}_S[b])_{s_0}-$module \ $E^{p+1}_y\big/A(E^{p+1}_y)$ \ is \ $b-$complete, for each \ $y \in \pi^{-1}(s_0) \cap \mathcal{Y}$.
\end{cor}

\parag{proof} Remark first that \ $E^{p+1}_y\big/A(E^{p+1}_y)$ \ is an \ $(\mathcal{O}_S[[a]])_{s_0}-$module of finite type  using the noetherianity of this ring and the fact that  the \ $(\mathcal{O}_S[[a]])_{s_0}-$ module  \ $\big(\Xi_{\Lambda, S}^{(n)}\big)_{s_0} \otimes_{\C} H^p(F_{y},\C)$ \ is of finite type. Remark then that the \ $a-$filtration  and \ $b-$filtration on \  $\big(\Xi_{\Lambda, S}^{(n)}\big)_{s_0} \otimes_{\C} H^p(F_{y},\C)$ \ are equivalent. So a Cauchy sequence for the \ $b-$filtration in \ $E^{p+1}_y\big/A(E^{p+1}_y)$ \ gives a Cauchy sequence for the \ $a-$filtration in \ $\big(\Xi_{\Lambda, S}^{(n)}\big)_{s_0} \otimes_{\C} H^p(F_{y},\C)$ \ and then converges for the \ $a-$filtration to an element which is in \ $E^{p+1}_y\big/A(E^{p+1}_y)$. Then it is enough to prove that the convergence for the \ $a-$filtration in \ $E^{p+1}_y\big/A(E^{p+1}_y)$ \ implies the convergence for the \ $b-$filtration. But this is a simple consequence of the hypothesis \ $(H0)$ \ because we have, for each \ $p \geq 0$, locally on \ $\mathcal{Y}$, an integer \ $N$ \  such that
$$ a^N.E^p \subset b.E^p $$
as an easy consequence of the Nullstellensatz : as \ $ \mathcal{Z} = \{ d_{/S}f = 0 \} \subset \mathcal{Y} = \{f = 0\}$ \ the support of the coherent sheaf \ $K^p_{/S}\big/ I^p_{/S} $ \ is contained in \ $\{f = 0\}$. So there exists locally on \ $\mathcal{Y}$ an integer \ $N$ \ such that \ $a^N.E^p \subset b.E^p$.$\hfill \blacksquare$

\parag{Remark} A consequence of the previous corollary is that \ $A(E_{y}^{p})$ \ is stable by \ $b$, and so we have the equality \ $\tilde{A}(E^{p}) = A(E^{p})$ \ for each \ $p \geq 0$. Then in the situation of the proposition \ref{b-separation} we shall have\footnote{see the proof of the proposition \ref{b-separation} for the first equality.} \ $B(E^{p}) = \tilde{A}(E^{p}) = A(E^{p})$ \ for each \ $p \geq 0$. So as in the proof above we obtain also locally on \ $\mathcal{Y}$ \ an integer \ $N$ \ such that \ $b^{2N}.B(E^{p}) = \{0\}$ \ for each \ $p \geq 0$.

 \subsection{The  theorem.}

Our first important result is the extension to the \ $S-$relative case  of the construction given in [B.II] theorem 2.1.1.

\begin{thm}\label{construction}
 In our situation \ $(@)$,  under the hypothesis \ $(H0), (H1), (H2)$, the map \ $ u : = u^1\circ u^0$ \ induces  for each \ $p \geq 0$ \ an isomorphism of \ $\mathcal{A}^0_S-$modules which is also \ $\mathcal{O}_S[[a]]-$linear
 $$ \mathcal{H}(u)^p : E^p \longrightarrow \mathcal{E}^p.$$
 Moreover we have the following properties of the \ $\mathcal{A}_S-$modules \ $\mathcal{E}^p, p \geq 0$ :\\
 \begin{enumerate}[i)]
 \item There exists locally on \ $\mathcal{Y}$ \ an integer \ $N$ \ such that \ $a^{N}$ \  kills \ $A(\mathcal{E}^p)$  \ for each \ $p \geq 0$.
 \item We have \ $B(\mathcal{E}^p) =  A(\mathcal{E}^p)$ \ for each \ $p \geq 0$.
  \item There exists locally on \ $\mathcal{Y}$ \ an integer \ $N'$ \ such that \ $b^{N'}$ \ kills the \ $b-$torsion \ $B(\mathcal{E}^p)$ \ of \ $\mathcal{E}^p$ \ for each \ $p \geq 0$.
  \item \ $\cap_{m \geq 0} \ b^m.\mathcal{E}^p = \{0\}$ \ for each \ $p \geq 0$.  
  \item For each \ $y \in \mathcal{Y}$ \ the germ \ $\mathcal{E}_y^{p+1}\big/A(\mathcal{E}^{p+1}_y)$ \ may be embeded in some \ $\mathcal{A}_{S,s_0}-$module of asymptotic expansions \ $\Xi_{\Lambda, S, s_0}^{(n)} \otimes V $ \ where \ $s_0  : = \pi(y)$, $V : = H^p(F_y, \C)$ \ and  $\Lambda \subset ]0,1]\cap \mathbb{Q}$ \ is a finite subset. Here \ $F_y$ \ is the Milnor fiber of \ $f$ \ at \ $y$, which coincides with the Milnor fiber at \ $y$ \  of \ $f$ \ restricted to \ $\pi^{-1}(s_0)$, thanks to the hypothesis \ $(H2)$. The set \ $\{ \exp(2i\pi.\lambda), \lambda \in \Lambda\}$ \ is the spectrum of the monodromy acting on \ $H^p(F_y, \C)$ .\\
\end{enumerate}
 \end{thm}

\parag{Remark} This last property \ $v)$ \ implies that the \ $\big(\mathcal{O}_S[[a]]\big)_{s_0}-$module  \ $E^p_y\big/A(E^p_y)$ \ (resp. \ $\big(\mathcal{O}_S[[b]])_{s_0}-$module  \ $E^p_y\big/A(E^p_y)$) is finitely generated by noetherianity (see the point  1. in  the begining of section 1.2.)$\hfill \square$\\

The proof of the theorem \ref{construction} will be  completed at the end of  the section 2.3.\\

We first remark that the cohomology sheaves \ $E^p$ \ satisfy the properties i) to v) for each \ $p \geq 0$ :\\

The property i) is proved in proposition  \ref{a-torsion} ; the properties  iv) is obtained in the proposition \ref{b-separation} and the property v) in the lemma  \ref{dev. asympt.}. The properties ii) and iii) are shown in the final remark of the section 2.2.\\

We shall now prove that the cohomology sheaves \ $\mathcal{E}^p$ \ satisfy the properties i) to v) under our hypotheses \ $(H0),(H1),(H2)$.

\subsubsection{The \ $b-$separation of  \ $\mathcal{E}^p$.}

We want to prove properties iii) and iv) for  the cohomology sheaves \ $\mathcal{E}^p$.

\begin{lemma}\label{b-torsion}
Assume that we have \ $b^{N}.B(E^p) = 0$ \ then we have \ $b^{N+1}.B(\mathcal{E}^p) = 0$.
\end{lemma}

\parag{proof} Let \ $X =  \sum_{j=0}^{\infty} \ b^j.x_j $ \ in \ $\mathcal{K}_{/S}^{p}$, $D_{/S}-$closed, such that \ $b^M.X = D_{/S}W$ \ with \ $M \geq N+1$.\\
 Define \ $U = \sum_{j=0}^{M-1} \ b^j.w_j$ \ and \ $V = \sum_{j=M}^{\infty} \ b^{j-M}.w_j$. Then we have
\begin{equation*}
 b^M.X = D_{/S}(U) + b^{M-1}.D_{/S}(b.V) \tag{@}
 \end{equation*}
and this implies \ $ D_{/S}(U) - b^{M-1}.(d_{/S}f\wedge w_M) =0 $ \ and so
$$ b^{M-1}.\big[b.X - D_{/S}(b.V)\big] = b^{M-1} .d_{/S}f\wedge w_M$$
and  \ $b.X = d_{/S}f\wedge w_M +  D_{/S}(b.V)$. Note that the relation

$$x_0 = -d_{/S}f\wedge w_{M-1} + d_{/S}(w_M),$$

consequence of \ $(@)$,   \  implies \ $d_{/S}f\wedge d_{/S}(w_M) = 0$, so we have  $$\mathcal{H}^p(u)[b.d_{/S}(w_M)] = [b.X].$$
 So \ $b^{M-1}.[d_{/S}(w_M)]$ \ is in the kernel of \ $\mathcal{H}^p(u)$ \ which is zero from our assumption. Then \ $b^N.[d_{/S}(w_M)] = 0$ \ and so \ $b^{N+1}.[X] = 0 $ \ in \ $\mathcal{E}^p$. $\hfill \blacksquare$\\

So the previous lemma and property iii) for the sheaves \ $E^{p}, \forall p \geq 0$ \ imply property iii) for the sheaves \ $\mathcal{E}^p, \forall p \geq 0$.\\

\begin{lemma}\label{b-separation2}
Assume now that we have, locally on \ $\mathcal{Y}$ \  an integer \ $N$ \ such that \ $b^N.B(\mathcal{E}^p) = 0$ \ then we have \ $\cap_{m \geq 0} b^m(\mathcal{E}^p) = \{0\}$.
\end{lemma}

\parag{Proof} Note first that the assertion is local on \ $\mathcal{Y}$.  Let \ $\Omega_0 \in \mathcal{K}_{/S}^p$ \ such that \ $D_{/S}\Omega_0 = 0 $ \ and such that for each \ $m \geq 0$ \ there exists \ $\Omega_m \in \mathcal{K}_{f/S}^p$ \ and \ $U_m \in \mathcal{K}_{/S}^{p-1}$ \ such that \ $D_{/S}\Omega_ m = 0 $ \ 
and  \ $b^m.\Omega_m = \Omega_0 + D_{/S}U_m $. Then we have for each \ $m \geq 0$
$$ b^m.(\Omega_m - b.\Omega_{m+1}) = D_{/S}(U_m - U_{m+1}) $$
and this means that \ $ [V_m] : = [\Omega_m] - b.[\Omega_{m+1}]$ \ is in \ $B(\mathcal{E}^p)$. As we know, thanks to our assumption, that there exists locally on \ $\mathcal{Y}$ \ an integer \ $N$ \ such that \ $b^{N}.B(\mathcal{E}^p) = \{0\}$, we obtain
\begin{equation*}
  \sum_{m =0}^{\infty} b^m.V_m = \sum_{m = 0}^{N-1} b^m.V_m + D_{/S}W  =   \sum_{m=0}^{\infty} ( b^m.\Omega_m - b^{m+1}.\Omega_{m+1})  = \Omega_0
\end{equation*}
where we have defined \ $W : = \sum_{m = N}^{\infty} \ b^m.W_m $ \ using the fact that for each \ $m \geq N$ \ we may write \ $b^m.V_m = D_{/S}W_m$ \ for some \ $W_m \in  \mathcal{K}_{/S}^{p-1}$, and where we use also that on a given open set \ $U \subset \mathcal{Y}$ \ the \ $\mathcal{O}_S(U)[[b]]-$module \ $\Gamma(U, \mathcal{K}_{/S}^{\bullet})$ \ is complete for its \ $b-$adic topology. So we find that \ $[\Omega_0]$ \ is in \ $B(\mathcal{E}^p)$. But then any \ $[\Omega_m]$ \ is in \ $B(\mathcal{E}^p)$ \ and for \ $m \geq N$ \ this gives \ $b^m.[\Omega_m] = [\Omega_0] = [0]$. $\hfill \blacksquare$.\\

So this second lemma, using the fact that property iii) is already obtained, gives the property iv) for the sheaves \ $\mathcal{E}^p, \forall p \geq 0$.\\

\subsubsection{The $b-$completion}

\begin{lemma}\label{Ker b}
Let \ $[X] \in \mathcal{E}^p$ \ such that \ $b.[X] = 0$. Then there exists \ $[x] \in E^p$ \ such that \ $b.[x] = 0$ \ and \ $\mathcal{H}^p(u)[x] = [X] $.
\end{lemma}

\parag{proof} Write \ $ X = \sum_{j=0}^{\infty} \ b^j.x_j $. Now if \ $b.X = D_{/S}U$ \ where \ $U : = \sum_{j=0}^{\infty} \ b^j.u_j$ \ is a local section of \ $\mathcal{K}^{p-1}_{/S}$ \ we shall have 
$$ 0 = d_{/S}(u_0) - d_{/S}f\wedge u_1\quad {\rm and } \quad u_0 \in K^{p-1}_{/S} .$$
So \ $d_{/S}f\wedge d_{/S}(u_1) = 0 $ \ and \ $d_{/S}(u_1)$ \ is in \ $K^p_{/S} \cap Ker\, d_{/S}$. Then we have in \ $E^p$ \ the relation \ $b.[d_{/S}(u_1)] = [d_{/S}(u_0)] = 0 $.\\
Write \ $U = u_0 + b.V$ \ and \ $V = u_1 + b.W$. Then \ $b.W$ \ is in \ $\mathcal{K}^{p-1}_{/S}$ \ and  we have
\begin{equation*}
 D_{/S}(U - u_0) = b.X - d_{/S}u_0 \quad {\rm and \ then} \quad X = D_{/S}V  
 \end{equation*}
But \ $V \not\in \mathcal{K}^{p-1}_{/S} $ \ in general. Nevertheless \ $ X - d_{/S}u_1 = D_{/S}(b.W) $ \
with \ $[d_{/S}(u_1)] $ \ is the image by \ $\mathcal{H}^p(u)$ \ of the class \ $d_{/S}(u_1)$ \ in \ $E^p$ \ which is in \ $ Ker \, b$. This conclude the proof. $\hfill \blacksquare$\\

\begin{cor}\label{B(E)}
For any integer \ $N \geq 1$ \ the map \ $\mathcal{H}^p(u)$ \ sends \ $Ker\, b^N$ \ in \ $E^p$ \ onto \ $Ker\, b^N$ \ in \ $\mathcal{E}^p$.
\end{cor}

\parag{proof} As the case \ $N=1$ \ is proved in the previous lemma it is enough to prove that if the statement is true for \ $N \geq 1$ \ (and for \ $N = 1$)  then it is true for \ $N+1$. Let \ $X \in Ker\, b^{N+1}$ \ in \ $\mathcal{E}^p$. Then \ $b.X$ \ may be written \ $\mathcal{H}^p(u)(x)$ \ where \ $x$ \ is in \ $E^p$ \ and satisfies \ $b^N.x = 0$. Now, as \ $\mathcal{H}^p(u)(x)$ \ is in \ $b.\mathcal{E}^p$ \ we may write \ $x = b.y  $ \  for some \ $y \in E^p$,  thanks to the lemma \ref{b-filt.} below. Then \ $Y : = \mathcal{H}^p(u)(y)$ \ satisfies \ $b.(X - Y) = 0$ \ and so we may find \ $t \in E^p$ \ such that \ $b.t = 0$ \ and \ $\mathcal{H}^p(u)(t) = X - Y $. Then \ $y + t$ \ is in \ $Ker\, b^{N+1}$ \ and is sent on \ $X$ \ by \ $\mathcal{H}^p(u)$. $\hfill \blacksquare$

\begin{lemma}\label{b-filt.}
Let \ $k$ \ be a positive integer. If \ $x \in E^p$ \ is such that  \ $\mathcal{H}(u)[x] $ \ is in \ $b^k.\mathcal{E}^p$ \ then \ $x$ \ is in \ $b^k.E^p$.
\end{lemma}

\parag{proof} We begin by the case \ $k = 1$. Then we may write
$$ x = b.X + D_{/S}U $$
with \ $X \in \mathcal{K}_{/S}^{p}, D_{/S}X = 0$ \ and \ $U \in \mathcal{K}_{/S}^{p-1}$ ; 
 this implies \ $x = d_{/S}(u_0) - d_{/S}f\wedge u_1$. So \ $[d_{/S}(u_1)]$ \ is in \ $E^p$ \ and \ $[d_{/S}(u_0)] = 0$ \ in \ $E^p$, as \ $u_{0} $ \ is in \ $K_{/S}^{p-1}$. This gives \ $[x] = -b[d_{/S}(u_1)]$, concluding this case.\\
Assume now \ $k \geq 1$ \ and that the case \ $k$ \ is proved. Then consider \ $[x] \in E^p$ \ such that
$$ x = b^{k+1}.X + D_{/S}U .$$
Then the inductive assumption gives \ $[y] \in E^p$ \ such that \ $[x] = b^k.[y]$. Then we have may assume directly  (up to change \ $U$) that
$$ b^k.y = b^{k+1}.X + .D_{/S}(U) . $$
This gives \ $d_{/S}u_{j} = d_{/S}f \wedge u_{j+1} \quad \forall j \in [0,k-1]$. So each \ $d_{/S}u_{j}$ \ for \ $j \in [0,k]$ \ defines an element in \ $E^{p}$, and we have (recall that \ $u_{0} \in K_{/S}^{p-1}$)
$$ 0 = [d_{/S}u_{0}] = b.[d_{/S}u_{1}] = \dots = b^{k}.[d_{/S}u_{k}] \quad {\rm in } \ E^{p} .$$
Now write \ $ U  = V + b^{k+1}.W$. The relation above becomes
$$ b^{k}.y = b^{k+1}.X + b^{k}.D_{/S}(b.W) + (d_{/S}u_{k}).b^{k} $$
which implies
$$ y - d_{/S}u_{k} = b.X + D_{/S}(b.W).$$
So we may apply the case \ $ k = 1$ \ to the class \ $[y - d_{/S}u_{k}] \in E^{p}$, and we find that \ $[y - d_{/S}u_{k}] \in b.E^{p}$. So \ $b^{k}.[y - d_{/S}u_{k}] = [x] + b^{k}.[d_{/S}u_{k}] $ \ is in \ $b^{k+1}.E^{p}$. Now the fact that \ $b^{k}.[d_{/S}u_{k}] = 0$ \ in \ $E^{p}$ \ allows to conclude.$\hfill \blacksquare$\\

\parag{Remark} This lemma implies that the \ $b-$filtration on \ $E^p$ \ is induced by the \ $b-$filtration on \ $\mathcal{E}^p$. $\hfill \square$

\subsubsection{The end of the proof.}

As we have proved the properties i) to v) for the sheaves \ $E^p$ \ the last point to prove is that the \ $\mathcal{A}^0_S-$linear map \ $\mathcal{H}^p(u) : E^p \to \mathcal{E}^p$ \ is bijective for each \ $p \geq 0$. But we already prove its injectivity thanks to the proposition  \ref{u1} using here condition iv) for the sheaves \ $E^{p}$ \ and also for the sheaf \ $\mathcal{E}^{p}$ \ thanks to  lemma \ref{b-separation2}.\\
Now the surjectivity comes from the \ $b-$density of its image proved in the proposition \ref{u1} and the following facts :
\begin{enumerate}
\item The map \ $\mathcal{H}^p(u)$ \ induce a bijective map between \ $B(E^p)$ \ and \ $B(\mathcal{E}^p)$.
\item The \ $b-$filtration on \ $E^p$ \ is induced by the \ $b-$filtration of \ $\mathcal{E}^p$.
\item For each \ $y \in \mathcal{Y}$ \ the quotient \ $\big(E^p\big/B(E^p)\big)_y$ \ is \ $b-$complete.
\item For each \ $p \geq 0$ \ $\mathcal{E}^p$ \ is \ $b-$separated (and so is \ $b-$complete).
\end{enumerate}
The point 1) is proved in the corollary \ref{B(E)}, thanks to lemma \ref{b-torsion} and the fact that there exists locally on \ $\mathcal{Y}$ \ an integer \ $N'$ \ such that \ $b^{N'}.E^{p} = \{0\}$.\\
The point 2) is proved in the lemma \ref{b-filt.}.\\
The point 3) is proved in the corollary \ref{completion 1}.   \\
The point 4) is proved in the lemmata  \ref{b-torsion}  and \ref{b-separation2}  thanks to the property iii)  for \ $E^p$.
This conclude the proof of the theorem \ref{construction}. $\hfill \blacksquare$

\subsubsection{Functorial properties.}

It is an easy exercice to prove the following functorial properties of the complex \ $ (\mathcal{K}^{\bullet}_{f/S}, D_{f/S}^{\bullet}) $ \ constructed above. 

\begin{prop}\label{functorial prop.}
In our situation \ $(@)$ \ with the hypothesis \ $(H0)$ \ the complex of \ $\mathcal{A}_S-$modules constructed in the theorem \ref{construction} has the following functorial properties:
\begin{enumerate}
  \item For any holomorphic map \ $\varphi : T \to S$ \ from a reduced complex space \ $T$ \ there exists a natural pull-back map of complexes over the pull-back of sheaves of algebras \ $\varphi^* : \mathcal{A}_S \to \mathcal{A}_T$
  $$ \varphi^* : (\mathcal{K}^{\bullet}_{f/S}, D^{\bullet}_{f/S}) \to (\mathcal{K}^{\bullet}_{g/T}, D^{\bullet}_{g/T}) $$
  where \ $g : \mathcal{X}\times_S T \to D$ \ is the composition of the first projection of the fiber product with \ $f$.  If \ $\varphi$ \ is a proper embedding, then \ $\varphi^*$ \ is given by the natural map to  the tensor product on \ $\mathcal{O}_S$ \ by \ $\mathcal{O}_T$.
  \item For any \ $S-$holomorphic map \ $\psi : \mathcal{X}_1 \to \mathcal{X}_2 $ \ between two \ $S-$manifolds and any holomorphic function \ $f : \mathcal{X}_2 \to D$ \ satisfying \ $(H0)$ \ and  such that \ $f\circ \psi$ \ satisfies \ $(H0)$, we have a natural pull-back of \ $\mathcal{A}_S-$complexes
  $$\psi^* : (\mathcal{K}^{\bullet}_{f/S}, D_{f/S}^{\bullet}) \to  (\mathcal{K}^{\bullet}_{g/S}, D_{g/S}^{\bullet})$$ 
  where \ $g : = f\circ \psi$.\\
  
  These two "pull-back" maps are of course compatible with the the obvious pull-back of relative differential forms and so with the natural map
  $$  u : (K_{/S}^{\bullet}, d_{/S}^{\bullet}) \rightarrow (\mathcal{K}_{/S}^{\bullet}, D_{/S}^{\bullet})    $$
  defined in the theorem \ref{construction}.
  \end{enumerate}
  \end{prop}

  \section{The finiteness theorem.}
  
   We assume in this section that we are in the situation \ $(@)$ \ with the hypotheses \ $(H0),(H1),(H2)$.

  \begin{defn}\label{petit}
  A left \ $\mathcal{A}_S-$module \ $E$ \ on \ $S$ \ is called {\bf $S-$small} when the following conditions hold
  \begin{enumerate}[1)]
  \item \ $\cap_{m\geq 0}\  b^m.E \subset A(E) $.
  \item \ $B(E) \subset A(E)$.
  \item Locally on \ $S$ \ there exists an integer \ $ N$ \ such that \ $ a^N.A(E) = 0 $.
  \item The sheaves \ $Ker\, b$ \ and \ $Coker\, b$ \ are \ $\mathcal{O}_S-$coherent.
  \end{enumerate}
  \end{defn}
  
  Remark that, using the results in section 2.2 \ the conditions 1) 2) and  3) implies that we have in fact \ $B(E) = \tilde{A}(E)$ \ and \ $b^{2N}.B(E) = \{0\}$ \ and then \ $\cap_{m\geq 0}\  b^m.E = \{0\}$.
  
  \begin{lemma}
  If a left \ $\mathcal{A}_S-$module \ $E$ \ is \ $S-$small then \ $E$ \ is a finite type \ $\mathcal{O}_S[[b]]-$module such that its \ $b-$torsion \ $B(E)$ \ is \ $\mathcal{O}_S-$coherent. 
  \end{lemma}
  
  \parag{proof} 
   We want to show that the sheaves \ $Ker\, b^{\nu}$ \ are \ $\mathcal{O}_S-$coherent for all \ $\nu \geq 1$. As this is true  for \ $\nu = 1$ \  by the condition 4), assume that it is proved for \ $\nu \geq 1$ \ and let us prove that \ $Ker\, b^{\nu+1}$ \ is \ $\mathcal{O}_S-$coherent. Consider the obvious \ $\mathcal{O}_S-$linear map \ $Ker\,b^{\nu} \to E\big/b.E$, between two \ $\mathcal{O}_S-$coherent sheaves. Its kernel is coherent, and as we have \ $Ker\, b^{\nu}\cap b.E = b(Ker\,b^{\nu+1})$ \ we obtain the coherence of
  \ $b(Ker\,b^{\nu+1})$. But now the exact sequence of \ $\mathcal{O}_S-$modules
  $$ 0 \to Ker\, b \to Ker\, b^{\nu+1} \overset{b}{\longrightarrow}  b(Ker\,b^{\nu+1}) \to 0 $$
  gives the coherence of \ $Ker\, b^{\nu+1}$. This implies the coherence of \ $B(E)$ \ thanks to the remark preceeding this lemma.\\
 Define now \ $F : = E\big/B(E) $. We shall prove the coherence of \ $F\big/b.F$. \\
  We have an isomorphism \ $F\big/b.F \simeq E\big/B(E) + b.E $ \ and an exact sequence
  $$0 \to B(E)\big/B(E) \cap b.E \to E\big/b.E \to F\big/b.F \to 0 .$$
  Now the equality  \ $B(E) \cap b.E = b.B(E)$ \ and the isomorphism  \ $b.B(E) \simeq B(E)\big/Ker\, b$ \ imply the coherence of \ $ B(E)\big/B(E) \cap b.E$ \ and the coherence of \ $E\big/b.E$ \ which is assumed, gives the coherence of \ $F\big/b.F$. \\
   But now the condition 1) implies that \ $F$ \ is \ $b-$complete and without \ $b-$torsion. So it is now standard to prove that \ $F$ \ is a  finite type \ $\mathcal{O}_S[[b]]-$module. We conclude that \ $E$ \ is  a finite type \ $\mathcal{O}_S[[b]]-$module.   $\hfill \blacksquare$
  
  \parag{Remark} We have proved that \ $F\big/b.F$ \ is also \ $\mathcal{O}_S-$coherent. Then there exists a Zariski dense open set \ $S' \subset S$ \ where \ $F\big/b.F$ \ is locally free, and, as \ $F$ \ has no \ $b-$torsion and is \ $b-$complete, it will be locally free over \ $\mathcal{O}_S[[b]]$ \ on \ $S'$. On such an open set, \ $E$ \ is locally a direct sum of its \ $b-$torsion with a free finite type \ $\mathcal{O}_S[[b]]-$module. $\hfill \square$\\

  \begin{defn}\label{geometric}
  We shall say that a left \ $\A-$module \ $E$ \ is {\bf geometric} when \ $E$ \ is \ $S-$small and when it associated family of (a,b)-module \ $E/B(E)$ \ has geometric fibers at each point of \ $S$.
  \end{defn}
  
  This definition is given in section 5.1 of  [B.II]. Of course, the formal completion of the Brieskorn module of a germ with an isolated singularity is always geometric. The condition to be geometric  encodes in the (a,b)-module setting the following general results on the Gauss-Manin connection of a holomorphic function : regularity, the monodromy theorem and Malgrange's positivity theorem.\\
  
  Let define the complex of \ $\A-$modules  \ $(\tilde{\mathcal{K}}_{/S}^{\bullet}, D_{/S}^{\bullet})$ \ on \ $\mathcal{Y}$ \ as the quotient of \ $(\mathcal{K}_{/S}^{\bullet}, D_{/S}^{\bullet})$ \ by the image by \ $u$ \ of the sub-complex (with \ $0$ \ differential)
 $$f^{-1}(\hat{\Omega}_{S\times D/S}^{1}) \subset (K_{/S}^{\bullet}, d_{/S}^{\bullet}).$$

  Our main result  is the following finiteness theorem, which shows that the Gauss-Manin connection produces in our situation a \ $\A-$module which is geometric .

 \begin{thm}\label{Finitude}
 Let \ $S$ \ be a reduced complex space and let \ $\mathcal{X}$ \ be a \ $S-$relative manifold of pure relative dimension \ $n + 1$. Note \ $\pi : \mathcal{X} \to S$ \ the projection. Let \ $f : =   \mathcal{X} \to  D$ \ be a   holomorphic function  to an open  disc \ $D$ \ in \ $\mathbb{C}$ \ with center \ $0$. Assume   that the properties  \ $(H0), (H1), (H2)$ \ are satisfied and that \ $\mathcal{Y}$ \ is \ $S-$proper. Then the complexes of \ $\A-$modules 
 $$\mathbb{R}\pi_*(\mathcal{K}_{f/S}^{\bullet}, D_{f/S}^{\bullet}) \quad {\rm and} \quad \mathbb{R}\pi_*(\tilde{\mathcal{K}}_{f/S}^{\bullet}, D_{f/S}^{\bullet})$$
 have   cohomology sheaves  which are geometric \ $\A-$modules for any degree.
 \end{thm}
 
 \parag{proof} Thanks to the theorem \ref{construction}, we may replace the complex \ $(\mathcal{K}_{f/S}^{\bullet}, D_{f/S}^{\bullet})$ \ by the complex \ $(K_{/S}^{\bullet}, d_{/S}^{\bullet})$. To show that \ $E : =  \mathcal{H}^p\big[\mathbb{R}\pi_*(K_{/S}^{\bullet}, d_{/S}^{\bullet})\big]$ \ is small, it is enough to prove that \ $E$ \ satisfies the condition 4) of the definition \ref{petit}, as the properties 1) 2) and 3) are given by the theorem \ref{construction}.\\
 Consider now the long exact sequence of hypercohomology of the exact sequence of complexes
 $$   0 \to (I_{/S}^{\bullet}, d_{/S}^{\bullet}) \to (K_{/S}^{\bullet}, d_{/S}^{\bullet}) \to ([K_{/S}/I_{/S}]^{\bullet}, d_{/S}^{\bullet}) \to 0 .$$
 It contains the exact sequence
 \begin{align*}
  & \dots \to  \mathcal{H}^{p-1}\big[\mathbb{R}\pi_*([K_{/S}\big/I_{/S}]^{\bullet}, d_{/S}^{\bullet})\big] \to \mathcal{H}^p\big[\mathbb{R}\pi_*(I_{/S}^{\bullet}, d_{/S}^{\bullet})\big] \overset{R^p(i)}{\to} \\
  & \quad \to  \mathcal{H}^p\big[\mathbb{R}\pi_*(K_{/S}^{\bullet}, d_{/S}^{\bullet})\big] \to \mathcal{H}^{p}\big[\mathbb{R}\pi_*([K_{/S}\big/I_{/S}]^{\bullet}, d_{/S}^{\bullet})\big] \to \dots
  \end{align*}
 and we know that \ $b$ \ is induced on the complex of \ $\A-$modules  quasi-isomorphic to \ $(K_{/S}^{\bullet}, d_{/S}^{\bullet})$ \ by the composition \ $i\circ \tilde{b}$ \ where \ $\tilde{b}$ \ is a quasi-isomorphism of complexes of \ $\mathcal{O}_S[[b]]-$modules (see lemma \ref{Tilde b} and the definition \ref{b en cohomologie}). This implies that the kernel  and the cokernel of \ $R^p(i)$ \ are isomorphic (as \ $\mathcal{O}_S-$modules) to \ $Ker\, b$ \ and \ $Coker \, b$ \ respectively. Now to prove that \ $E$ \ satisfies condition 4) of the definition \ref{petit} it is enough to prove coherence of the \ $\mathcal{O}_S-$modules  \ $ \mathcal{H}^{j}\big[\mathbb{R}\pi_*([K_{/S}\big/I_{/S}]^{\bullet}, d_{/S}^{\bullet})\big] $ \ for all \ $j \geq 0 $.\\
 But the sheaves \ $[K_{/S}\big/I_{/S}]^j \simeq [Ker\,d_{/S}f\big/Im\, d_{/S}f]^j$ \ are coherent on \ $\mathcal{X}$ \ and supported in \ $\mathcal{Y}$. The spectral sequence
 $$ E_2^{p,q} : = H^q\big( R^p\pi_*( [K_{/S}\big/I_{/S}]^{\bullet}), d_{/S}^{\bullet}\big) $$ 
 which converges to \ $ \mathcal{H}^{j}\big[\mathbb{R}\pi_*([K_{/S}\big/I_{/S}]^{\bullet}, d_{/S}^{\bullet})\big] $, is a bounded complex of coherent  \ $\mathcal{O}_S-$modules by the direct image theorem of H. Grauert. This gives the desired finiteness.\\
 To conclude the proof, we want to show that \ $E/B(E)$ \ is geometric. But this is an easy consequence of the regularity of the Gauss-Manin connexion of \ $f$, the monodromy theorem and Malgrange's positivity theorem. \\
 The case of the complex \ $(\tilde{\mathcal{K}}_{/S}^{\bullet}, D_{/S}^{\bullet})$ \ follows in the same way.$\hfill \blacksquare$\\

\parag{Remark} 
  As the complex \ $(\tilde{\mathcal{K}}_{/S}^{\bullet}, D_{/S}^{\bullet})$ \  has its cohomology supported in \ $\mathcal{Z}$ \ under the hypothesis \ $(H0)$, to obtain that  \ $\mathbb{R}\pi_*(\tilde{\mathcal{K}}_{f/S}^{\bullet}, D_{f/S}^{\bullet})$ \  are geometric \ $\A-$modules for any degree, under the hypotheses \ $(H0),(H1),(H2)$, the assumption that the restriction of \ $f$ \ to \ $\mathcal{Z}$ \ is proper would be enough.$\hfill \square$

  \section{Bibliography}

\begin{itemize}

\item{[Br.70]} Brieskorn, E. {\it Die Monodromie der Isolierten Singularit{\"a}ten von Hyperfl{\"a}chen}, Manuscripta Math. 2 (1970), p. 103-161.


\item{[B.93]} Barlet, D. {\it Th\'eorie des (a,b)-modules I}, in Complex Analysis and Geo-metry, Plenum Press, (1993), p. 1-43.

\item{[B.95]} Barlet, D. {\it Th\'eorie des (a,b)-modules II. Extensions}, in Complex Analysis and Geometry, Pitman Research Notes in Mathematics Series 366 Longman (1997), p. 19-59.


\item{[B.I]} Barlet, D. {\it Sur certaines singularit\'es non isol\'ees d'hypersurfaces I}, Bull. Soc. math. France 134 (2), ( 2006), p.173-200.

\item{[B.II]} Barlet, D. {\it Sur certaines singularit\'es d'hypersurfaces II}, J. Alg. Geom. 17 (2008), p. 199-254.





\item{[B.-S. 04]} Barlet, D. et Saito, M. {\it Brieskorn modules and Gauss-Manin systems for non isolated hypersurface singularities,} J. Lond. Math. Soc. (2) 76 (2007) $n^01$ \ p. 211-224.


\item{[D.70]} Deligne, P. {\it Equations diff\'erentielles \`a points singuliers r\'eguliers} Lecture Notes in Math. 163, Springer 1970.

\item{[G.65]} Grothendieck, A. {\it On the de Rham cohomology of algebraic varieties} Publ. Math. IHES 29 (1966), p. 93-101.


\item{[M.74]} Malgrange, B. {\it Int\'egrale asymptotique et monodromie}, Ann. Sc. Ec. Norm. Sup. 7 (1974), p. 405-430.


\item{[Mi.68]  Milnor, J.  \textit{Singular Points of Complex Hypersurfaces .} Ann. of Math. Studies  61 (1968) Princeton  .}

\item{[S.89]} Saito, M. {\it On the structure of Brieskorn lattices}, Ann. Inst. Fourier 39 (1989), p. 27-72.

\end{itemize}

\end{document}